\documentclass[10pt,a4paper]{amsart}
\date{June 9, 2011}      
\usepackage{hyperref}
\usepackage{latexsym,amsmath,amsfonts,amscd,amssymb,graphics}
\usepackage[all]{xy}
\usepackage[ansinew]{inputenc}

\theoremstyle{plain}  
\newtheorem{theorem}{Theorem}[section]
\newtheorem*{theorem*}{Theorem}

\newtheorem{corollary}[theorem]{Corollary}
\newtheorem{lemma}[theorem]{Lemma}
\newtheorem{proposition}[theorem]{Proposition}

\theoremstyle{definition}
\newtheorem{definition}[theorem]{Definition}
\newtheorem{notation}[theorem]{Notation}

\theoremstyle{remark}

\newtheorem{remark}[theorem]{Remark}

\newtheorem*{claim*}{Claim}

\numberwithin{equation}{section}

\newcommand{\suchthat}{\;|\;}
\newcommand{\st}{\;|\;}
\newcommand{\abs}[1]{\lvert#1\rvert}
\renewcommand{\leq}{\leqslant}

\renewcommand{\geq}{\geqslant}

\renewcommand{\setminus}{\smallsetminus}
\newcommand{\x}{\times}
\newcommand{\into}{\hookrightarrow}

\newcommand{\R}{\mathbb{R}}

\newcommand{\Z}{\mathbb{Z}}
\newcommand{\C}{\mathbb{C}}

\newcommand{\cM}{\mathcal{M}}
\newcommand{\cN}{\mathcal{N}}

\newcommand{\cR}{\mathcal{R}}
\newcommand{\cS}{\mathcal{S}}

\newcommand{\lie}{\mathfrak}

\newcommand{\PGL}{\mathrm{PGL}}

\newcommand{\Sp}{\mathrm{Sp}}
\newcommand{\PSp}{\mathrm{PSp}}

\newcommand{\U}{\mathrm{U}}
\newcommand{\EGL}{\mathrm{EGL}}
\newcommand{\ESp}{\mathrm{ESp}}

\newcommand{\EU}{\mathrm{EU}}

\newcommand{\GL}{\mathrm{GL}}

\newcommand{\cSp}{\mathrm{Sp}}

\DeclareMathOperator{\Jac}{Jac}

\DeclareMathOperator{\divisor}{div}

\DeclareMathOperator{\cSym}{Sym}
\DeclareMathOperator{\cSupp}{Supp}
\DeclareMathOperator{\ord}{ord} 
 
\DeclareMathOperator{\tr}{tr}
\DeclareMathOperator{\rk}{rk} 

\DeclareMathOperator{\im}{im} 

\DeclareMathOperator{\Hom}{Hom} 
\DeclareMathOperator{\End}{End}
\DeclareMathOperator{\Ext}{Ext}

\DeclareMathOperator{\codim}{codim}

\newcommand{\Pic}{\operatorname{Pic}}

\newcommand{\Aut}{\operatorname{Aut}}

\newcommand{\liem}{\mathfrak{m}}
\newcommand{\liez}{\mathfrak{z}}

\newcommand{\liemc}{\mathfrak{m}^{\mathbb{C}}}
\newcommand{\lieh}{\mathfrak{h}}
\newcommand{\liehc}{\mathfrak{h}^{\mathbb{C}}}
\newcommand{\lieg}{\mathfrak{g}}
\newcommand{\liegc}{\mathfrak{g}^{\mathbb{C}}}

\let\oldmarginpar\marginpar
\renewcommand\marginpar[1]{\oldmarginpar{\tiny\bf\begin{flushleft} #1
\end{flushleft}}}

\begin{document}

\title[Quadratic pairs and surface group
  representations]{Rank two quadratic pairs and surface group
  representations}

\author[P. B. Gothen]{Peter B. Gothen}
\address{Centro de Matem\'atica da Universidade do Porto\\
Faculdade de Ci\^encias, Universidade do Porto \\
Rua do Campo Alegre \\ 4169-007 Porto \\ Portugal }
\email{pbgothen@fc.up.pt}

\author[A. G. Oliveira]{Andr\'e G. Oliveira}
\address{Centro de Matem\'atica\\
  Universidade de Tr\'as-os-Montes e Alto Douro\\
  Quinta dos Prados, Apartado 1013 \\ 5000-911 Vila Real \\ Portugal }
\email{agoliv@utad.pt}

\thanks{
Members of VBAC (Vector Bundles on Algebraic Curves).
Partially supported by CRUP through Acção Integrada Luso-Espanhola nº
E-38/09
and by FCT (Portugal) through the projects PTDC/MAT/099275/2008 and
PTDC/MAT/098770/2008, and through the Centro de Matemática da
Universidade do Porto (first author) and Centro de Matemática da
Universidade de Trás-os-Montes e Alto Douro (second author).
}


\begin{abstract}
  Let $X$ be a compact Riemann surface. A quadratic pair on $X$
  consists of a holomorphic vector bundle with a quadratic form which
  takes values in fixed line bundle. We show that the moduli spaces of
  quadratic pairs of rank $2$ are connected under some constraints on
  their topological invariants. As an application of our results we
  determine the connected components of the $\mathrm{SO}_0(2,3)$-character
  variety of $X$.
\end{abstract}

\maketitle


\section{Introduction}

Let $X$ be a compact Riemann surface of genus $g\geq 2$. Many kinds of pairs
$(V,\varphi)$ on $X$, consisting of a holomorphic vector bundle $V \to
X$ and a holomorphic section $\varphi$ of an associated bundle, have
been extensively studied. Important examples are Bradlow pairs \cite{bradlow:1991}, where $\varphi \in H^0(X,E)$ lives in the fundamental
representation and Higgs bundles \cite{hitchin:1987}, where $\varphi
\in H^0(X,\End(E) \otimes K)$ lives in the adjoint representation
(twisted by the canonical bundle $K$ of $X$).  Many
more examples can be found in the survey \cite{bradlow-et-al:1995}.

In this paper we focus on \emph{$U$-quadratic pairs} $(V,\gamma)$,
where $\gamma$ is a global section of $S^2V^*\otimes U$ for a fixed
line bundle $U \to X$. These are of interest for at least
two reasons. On the one hand they can be viewed as giving rise to
bundles of  quadrics and hence form a very natural
generalization of the linear objects of vector bundles. On the other
hand they arise naturally in the study of another kind of linear pairs,
namely $G$-Higgs bundles: these are the appropriate objects for
studying character varieties for the fundamental group of $X$ in a
real Lie group $G$ through the non-abelian Hodge theory correspondence
(see for example \cite{bradlow-garcia-prada-gothen:2005} for a survey
on this topic).

Moduli spaces of quadratic pairs were constructed via GIT and studied
by Gómez and Sols in \cite{gomez-sols:2000} and also by Schmitt in
\cite{schmitt:2004}. Moreover, Mundet in the appendix to
\cite{gomez-sols:2000}, showed that the stability condition used for
constructing moduli is the same one which allows to prove a
Hitchin--Kobayashi correspondence for quadratic pairs, relating
stability of the quadratic pair to the existence of solutions to
certain gauge theoretic equations.  This stability condition depends
on a real parameter $\alpha$ hence, for each value of this parameter,
there is a moduli space which we denote by $\cN_\alpha(n,d)$.

In the first part of this paper, we study the number of connected
components of the moduli spaces of $U$-quadratic pairs on $X$. Our
strategy is the one pioneered by Thaddeus \cite{thaddeus:1994} and
subsequently used in many other cases, e.g.\
\cite{bradlow-garcia-prada-gothen:2004 triples}. It consists in
studying the variation of the moduli space $\cN_\alpha(n,d)$ with the
parameter $\alpha$. As usually happens, when we run over $\alpha$, the
moduli spaces $\cN_\alpha(n,d)$ are isomorphic for parameter values in
intervals and only change at a discrete set of \emph{critical
  values}. In these cases, the difference between the moduli spaces
are confined to subvarieties, which are called the \emph{flip
  loci}. For $n=2$, we describe explicitly these subvarieties and show
that they have positive codimension in $\cN_\alpha(2,d)$.  A necessary
condition for the non-emptiness of $\cN_\alpha(2,d)$ is $\alpha\leq
d/2$. Moreover, if $d>d_U=\deg(U)$, then $\cN_\alpha(2,d)$ is empty
unless $\alpha=d/2$ and $\cN_{d/2}(2,d)$ is the moduli space of
semistable rank $2$ and degree $d$ vector bundles. So we consider only
$d<d_U$ (the $d=d_U$ case is special and not considered here).  We
show that there is an $\alpha_m$ such that the $\cN_\alpha(2,d)$'s for
$\alpha<\alpha_m$ are all isomorphic. Then, using the theory of the
Hitchin system, and in particular the results obtained in
\cite{gothen-oliveira:2010}, we show that $\cN_{\alpha_m^-}(2,d)$ is
connected, where $\alpha_m^-$ is any value less than $\alpha_m$. This,
together with study of the flip loci, provides a similar conclusion
for the connectedness of the other spaces $\cN_\alpha(2,d)$, whenever
$d_U-d>g-1$ holds. Our result (Theorem
\ref{Nalpha(d,2) connected}) states then the following:

\begin{theorem*}
  Let $d$ and $d_U$ be such that $d_U-d>g-1$. Then, for every
  $\alpha\leq d/2$, the moduli space $\cN_\alpha(2,d)$ is connected.
\end{theorem*}

In the last part of the paper we apply our results to counting the
connected components of the character variety
\begin{displaymath}
  \cR(\pi_1X,\mathrm{SO}_0(2,3))
  =\Hom^{\mathrm{red}}(\pi_1X,\mathrm{SO}_0(2,3))/\mathrm{SO}_0(2,3),
\end{displaymath}
i.e., the space of reductive representations $\rho\colon \pi_1X \to
\mathrm{SO}_0(2,3)$ modulo the action by simultaneous
conjugation. Such a representation $\rho$ has a topological invariant
$c(\rho)=(\tau(\rho),w_2(\rho)) \in \pi_1\mathrm{SO}_0(2,3) \cong \Z
\x \Z/2$ given by the topological class of the associated flat
bundle. Now, for representations of $\pi_1 X$ into any isometry group
of a hermitean symmetric space of non-compact type there is an integer
invariant, known as the \emph{Toledo invariant}, and in the present
case this invariant is just the first coordinate $\tau(\rho)$ of
the topological class. Moreover, the Toledo invariant is bounded by the
Milnor--Wood type inequality
\begin{displaymath}
  \abs{\tau(\rho)} \leq 2g-2.
\end{displaymath}
For $(a,w) \in \Z \x \Z/2$, denote by
$\cR_{a,w}(\pi_1X,\mathrm{SO}_0(2,3)) \subseteq
\cR(\pi_1X,\mathrm{SO}_0(2,3))$ the subspace of representations $\rho$
such that $c(\rho)=(a,w)$.

Our results on quadratic pairs then lead to the following (Theorem
\ref{connected-components-RSO23}).

\begin{theorem*}
  For each $(a,w)\in\Z\times\Z/2$ such that $0<\abs{a}< 2g-2$, the space
  $\cR_{a,w}(\pi_1X,\mathrm{SO}_0(2,3))$ is connected.
\end{theorem*}

To put this result into perspective, recall that $\mathrm{SO}_0(2,3)$
is isomorphic to the adjoint form $\PSp(4,\R)$ of the real symplectic
group $\Sp(4,\R)$. It can be seen that a representation $\rho\colon
\pi_1 X \to \mathrm{SO}_0(2,3)$ lifts to $\Sp(4,\R)$ if and only if
$\tau(\rho) = w_2(\rho) \mod 2$. Moreover, if this is the case,
$\tau(\rho)$ coincides with the Toledo invariant of the lifted
representation and this in turn coincides with the topological
invariant in $\pi_1\Sp(4,\R) \cong \Z$ of
the lifted representation. The
connected components of the character variety
$\cR_{a}(\pi_1X,\Sp(4,\R))$ of representations in $\Sp(4,\R)$ with
topological invariant $a \in \Z$ have been determined in
\cite{gothen:2001} for $\abs{a}=0$ and $\abs{a}=2g-2$, and for the
remaining values of $\abs{a}$ by Garc\'\i{}a--Prada and Mundet in
\cite{garcia-prada-mundet:2004}. In the case of representations which
lift to $\Sp(4,\R)$ these results easily lead to the count of
connected components for representations in $\mathrm{SO}_0(2,3)$ (cf.\
\cite{bradlow-garcia-prada-gothen:2005}). Thus our Theorem completes
the count of the connected components for representations in
$\mathrm{SO}_0(2,3)$ and the novelty lies in the cases $ a \neq w
\mod 2$.

This paper is organized as follows. In
Section~\ref{sec:quadratic-pairs} we recall some basic facts about
quadratic pairs. In Section \ref{sec:variation-of-moduli} we carry out
the analysis of the variation of the moduli spaces $\cN_\alpha(2,d)$
with the parameter, leaving however the proof of connectedness of
$\cN_{\alpha_m^-}(2,d)$ for Section~\ref{The space
  Nalpham-(2,d)}. Then, in
Section~\ref{Number of connected components}, we put our results
together to obtain the main connectedness theorem for the moduli of
quadratic pairs. Finally, in Section~\ref{sec:higgs-surface-groups},
we give the application of our results to Higgs bundles and
representations of surface groups in the group $\mathrm{SO}_0(2,3)$.

\section{Quadratic pairs}
\label{sec:quadratic-pairs}
\subsection{Quadratic pairs and their moduli spaces}

Let $X$ be smooth projective curve over $\C$ of genus $g\geq 2$, and
let $U$ be a fixed holomorphic line bundle over
$X$. Write $$d_U=\deg(U)$$ for the degree of $U$.

\begin{definition}\label{definition of U quadratic pair}
  A \emph{$U$-quadratic pair} on $X$ is a pair $(V,\gamma)$, where $V$
  is a holomorphic vector bundle over $X$ and $\gamma$ is a global
  holomorphic non-zero section of $S^2V^*\otimes U$, i.e., $\gamma\in
  H^0(X,S^2V^*\otimes U)$.  The \emph{rank} and \emph{degree} of a
  quadratic pair are the rank and degree of the underlying vector
  bundle $V$, respectively. We say that $(V,\gamma)$ is of \emph{type}
  $(n,d)$ if $\rk(V)=n$ and $\deg(V)=d$.
\end{definition}

We shall often refer to a $U$-quadratic pair simply as a
\emph{quadratic pair}. Quadratic pairs are sometimes called
\emph{conic bundles} in the literature.

\begin{definition}
  Two $U$-quadratic pairs $(V,\gamma)$ and $(V',\gamma')$ are
  \emph{isomorphic} if there is an isomorphism $f:V\to V'$ such that
  $\gamma'f=((f^t)^{-1}\otimes 1_U)\gamma$, i.e., such that the
  following diagram commutes:
$$\xymatrix{V\ar[r]^{f}\ar[d]_{\gamma}&V'\ar[d]^{\gamma'}\\
V^*\otimes U\ar[r]_{(f^t)^{-1}\otimes 1_U}&V'^*\otimes U.}$$
\end{definition}

Quadratic pairs of rank $n \leq 3$ were studied in
\cite{gomez-sols:2000} by Gómez and Sols. They introduced an
appropriate $\alpha$-semistability condition, depending on a real
parameter $\alpha$\footnote{In fact a different parameter $\tau$ is
  used in \cite{gomez-sols:2000}. In the cases of interest to us, the
  precise definition of $\alpha$-semistability is given below, as well
  as the relation between the parameters $\tau$ and $\alpha$.}, and
constructed moduli spaces of $S$-equivalence classes of
$\alpha$-semistable quadratic pairs using GIT. The construction of the
moduli spaces for general rank is due to
Schmitt~\cite{schmitt:2004,schmitt:2008}. We denote the moduli space
of $S$-equivalence classes of $\alpha$-semistable $U$-quadratic pairs
on $X$ of rank $n$ and degree $d$ by
\begin{displaymath}
  \cN_{X,\alpha}(n,d)=\cN_\alpha(n,d).
\end{displaymath}

There is a Hitchin--Kobayashi correspondence for quadratic pairs. This
follows from the general results of \cite{mundet:2000},
\cite{bradlow-garcia-prada-mundet:2003} and
\cite{garcia-gothen-mundet:2008}; the Appendix to
\cite{gomez-sols:2000} treats the application to the case of quadratic
pairs. It says the a quadratic pair supports a solution to a certain
natural gauge theoretic equation if and only if it is
\emph{$\alpha$-polystable} (see below for the definition of this
concept). Moreover, each $S$-equivalence class has a unique
$\alpha$-polystable representative and thus $S$-equivalence of
$\alpha$-polystable pairs reduces to isomorphism. We can therefore
also consider $\cN_\alpha(n,d)$ as the moduli space of isomorphism
classes of $\alpha$-polystable quadratic pairs.

\subsection{\texorpdfstring{$U$}{U}-quadratic pairs of rank \texorpdfstring{$1$}{1}}

Although we will be mainly interested in quadratic pairs of type
$(2,d)$, we shall also need the description of the moduli spaces of
quadratic pairs on $X$ of rank $1$ and of their moduli spaces.

\begin{definition}\label{sspairs1}
  Fix a real parameter $\alpha$. A $U$-quadratic pair $(L,\delta)$ of
  type $(1,d)$ is \emph{$\alpha$-stable} if $\alpha\leq d$.
\end{definition}

This definition is equivalent to the one used in
\cite{gomez-sols:2000} of $\tau$-semistability for rank $1$ pairs. The
equivalence is obtained by taking $\tau=d-\alpha$. 

\begin{remark}\label{degenerated cases}
There are no strictly $\alpha$-semistable quadratic pairs of rank $1$.
\end{remark}

For quadratic pairs of type $(1,d')$, all the moduli spaces
$\cN_\alpha(1,d')$ with $\alpha\leq d'$ are isomorphic and there is
only one so-called \emph{critical value} of $\alpha$, for which the
stability condition changes, namely $\alpha=d'$.

\begin{lemma}
Let $\cN_\alpha(1,d')$ be the moduli space of $\alpha$-stable
quadratic pairs of type $(1,d')$. Then,
\begin{enumerate}
 \item for all $\alpha>d'$, $\cN_\alpha(1,d')=\emptyset$;
 \item for all $d'>d_U/2$, $\cN_\alpha(1,d')=\emptyset$.
\end{enumerate}
\end{lemma}
\proof The first item follows from the stability condition. For the
second part, we have that, if $(M,\delta)\in\cN_\alpha(1,d')$, then
$\delta:M\to M^{-1}U$ is non-zero and holomorphic, so $-2d'+d_U\geq
0$.
\endproof 

\begin{proposition}\label{description of N(1,dU/2)}
  Suppose that $d_U$ is even. If $\alpha\leq d_U/2$, then the moduli
  space $\cN_\alpha(1,d_U/2)$ is isomorphic to 
  $$
  S=\{F\in\Pic^{d_U/2}(X)\suchthat F^2\cong U\},
  $$
  the set of the $2^{2g}$ square roots of $U$.
\end{proposition}

\proof 
Let $\alpha\leq d_U/2$. If $(M,\delta)\in\cN_\alpha(1,d_U/2)$, then
$\delta:M\to M^{-1}U$ must be non-zero, hence an isomorphism.
Moreover, it is defined up to a non-zero scalar so the map
$\cN_\alpha(1,d_U/2)\to S$, $(M,\delta)\mapsto M$ is an isomorphism.
\endproof

It remains to describe $\cN_\alpha(1,d')$ for $\alpha\leq d'$ and
$d'<d_U/2$. Denote by $\cSym^n(X)$ the $n$th symmetric product of $X$,
the smooth variety which parametrizes the degree $n$ effective divisors on $X$.

\begin{proposition}\label{description of Nalpha(1,d) as a cover}
  If $d'<d_U/2$ and $\alpha\leq d'$, then $\cN_\alpha(1,d')$ is the
  $2^{2g}$-fold cover of the symmetric product $\cSym^{d_U-2d'}(X)$
  obtained by pulling back, via the Abel-Jacobi map, the cover
  $\Pic(X) \to \Pic(X)$ given by squaring of line bundles.
\end{proposition}

\proof Consider the map $\pi:\cN_\alpha(1,d')\to \cSym^{d_U-2d'}(X)$,
$(M,\delta)\mapsto\divisor(\delta)$.  Given $D\in\cSym^{d_U-2d'}(X)$,
$\pi^{-1}(D)$ is isomorphic to the set of square roots of $U(-D)$. The
result follows.
\endproof

\begin{corollary}\label{dimNalpha(1,d)}
 Let $d'<d_U/2$ and $\alpha\leq d'$. Then $\dim\cN_\alpha(1,d')=d_U-2d'$.
\end{corollary}

\subsection{Stability of quadratic pairs of rank $2$}

Our main objects of interest are type $(2,d)$ quadratic pairs
$(V,\gamma)$, where $V$ is a holomorphic vector bundle of rank $2$ and
degree $d$, and $\gamma\in H^0(X,S^2V^*\otimes U)$.  Most of
the time we will think of $\gamma$ as a holomorphic map
$\gamma:V\longrightarrow V^*\otimes U$ which is symmetric, i.e.,
$\gamma^t\otimes 1_U=\gamma$.

Given a rank $2$ vector bundle $V$ and a line subbundle $L\subset V$,
we denote by $L^\perp$ the kernel of the projection $V^*\to L^{-1}$. It
is a line subbundle of $V^*$ and $V/L$ is isomorphic to
$(L^\perp)^{-1}$.

The general definition of stability from
\cite{mundet:2000,bradlow-garcia-prada-mundet:2003,garcia-gothen-mundet:2008}
specializes as follows in the case of quadratic pairs. It is
equivalent to the definition of $\tau$-semistability of Gómez and
Sols~\cite{gomez-sols:2000} by taking $\tau=d/2-\alpha$.

\begin{definition}\label{sspairs2}
  Fix $\alpha\in\R$. A $U$-quadratic pair $(V,\gamma)$ of type
  $(2,d)$ is:
\begin{itemize}
\item \emph{$\alpha$-semistable} if $\alpha\leq d/2$ and, for every line
  subbundle $L\subset V$,
\begin{enumerate}
 \item $\deg(L)\leq\alpha$ if $\gamma\in H^0(X,(L^\perp)^2U)$;
 \item $\deg(L)\leq d/2$ if $\gamma\in H^0(X,L^\perp\otimes_S V^*\otimes U)$;
 \item $\deg(L)\leq d-\alpha$ if $\gamma\notin H^0(X,L^\perp\otimes_S V^*\otimes U)$.
\end{enumerate}
\item \emph{$\alpha$-stable} if it is $\alpha$-semistable and strict
  inequalities hold in $(1)$, $(2)$ and $(3)$ above.
\item \emph{$\alpha$-polystable} if $\alpha\leq d/2$ and, for every line
  subbundle $L\subset V$,
\begin{enumerate}
\item $\deg(L)\leq\alpha$ if $\gamma\in
  H^0(X,(L^\perp)^2U)$. Moreover, if $\deg(L)=\alpha$, there is
  $L'\subset V$ such that $V=L\oplus L'$;
\item $\deg(L)\leq d/2$ if $\gamma\in H^0(X,L^\perp\otimes_S
  V^*\otimes U)$. Moreover, if $\deg(L)=d/2$, there is $L'\subset V$
  such that $V=L\oplus L'$ and $\gamma'\in
  H^0(X,L^{-1}L'^{-1}U)$ such that $\gamma=\gamma'\oplus\gamma'$;
 \item $\deg(L)\leq d-\alpha$ if $\gamma\notin H^0(X,L^\perp\otimes_S
   V^*\otimes U)$. Moreover, if $\deg(L)=\alpha$, there is $L'\subset
   V$ such that $V=L\oplus L'$.
\end{enumerate}
\end{itemize}
\end{definition}

\begin{remark}\label{d/2-semistability}
The $d/2$-(semi)stability condition for $(V,\gamma)$ is equivalent to the usual
(semi)stability condition for the vector bundle $V$.
\end{remark}

An $\alpha$-semistable quadratic pair $(V,\gamma)$ is \emph{strictly
  $\alpha$-semistable} if it is not $\alpha$-stable. From the previous definition, we can separate strictly $\alpha$-semistable
quadratic pairs into three types.

\begin{definition}
  An $\alpha$-semistable quadratic pair $(V,\gamma)$ is strictly
  $\alpha$-semistable of \emph{type}:
  \begin{enumerate}
  \item[\textbf{(A)}] if there is a holomorphic line bundle $L\subset V$
    such that $\gamma\in H^0(X,(L^\perp)^2U)$ and $\deg(L)=\alpha$;
  \item[\textbf{(B)}] if there is a holomorphic line bundle $L\subset V$
    such that $\gamma\in H^0(X,L^\perp\otimes_S V^*\otimes U)$ and
    $\deg(L)=d/2$;
  \item[\textbf{(C)}] if there a holomorphic line bundle $L\subset V$
    such that $\deg(L)=d-\alpha$.
  \end{enumerate}
\end{definition}

\begin{definition}
  For a given type $(2,d)$, the values of $\alpha$ for which strictly
  $\alpha$-semistable quadratic pairs of type \textbf{(A)} or
  \textbf{(C)} exist are called \emph{critical values}, and the other
  values of $\alpha$ are called \emph{generic values}.
\end{definition}

\begin{remark}
  For generic $\alpha$ and for a pair $(V,\gamma)$, if there is no
  $L\subset V$ such that $\gamma(L)\subset L^\perp U$ and
  $\deg(L)=d/2$, then $(V,\gamma)$ is $\alpha$-semistable if and only
  if it is $\alpha$-stable. In particular, if $d$ is odd there are no strictly $\alpha$-semistable pairs of
  type \bf{(B)}.
\end{remark}

\begin{lemma}\label{sub}Let $(V,\gamma)$ be a $U$-quadratic pair of
  rank $2$ and let $L$ be a line subbundle of $V$. Then,
  \begin{enumerate}
  \item $\gamma\in
    H^0(X,(L^\perp)^2U)\Longleftrightarrow\gamma(L)=0\Longleftrightarrow\gamma(V)\subset
    L^\perp U$;
 \item  $\gamma\in H^0(X,L^\perp\otimes_SV^*\otimes U)\Longleftrightarrow\gamma(L)\subset L^\perp U$.
\end{enumerate}
\end{lemma}
\proof
This is an exercise in fibrewise linear algebra; see
\cite{oliveira:2008} for details.
%
%
\endproof

Using this lemma, we can rewrite the $\alpha$-(poly,semi)stability
condition in the following way.

\begin{proposition}\label{simplified notion of polystability}
Let $(V,\gamma)$ be a quadratic pair.
\begin{itemize}
\item The pair $(V,\gamma)$ is $\alpha$-semistable if and only if
  $\alpha\leq d/2$ and, for any line bundle $L\subset V$, the
  following conditions hold:
\begin{enumerate}
\item $\deg(L)\leq\alpha$, if $\gamma(L)=0$;
\item $\deg(L)\leq d/2$, if $\gamma(L)\subset L^\perp U$;
\item $\deg(L)\leq d-\alpha$, if $\gamma(L)\not\subset L^\perp U$.
\end{enumerate}
\item The pair $(V,\gamma)$ is $\alpha$-stable if and only if it is
  $\alpha$-semistable for any line bundle $L\subset V$, the conditions
  $(1)$, $(2)$ and $(3)$ above hold with strict inequalities.
\item The pair $(V,\gamma)$ is $\alpha$-polystable if and only if
  $\alpha\leq d/2$ and, for any line bundle $L\subset V$, the
  following conditions hold:
\begin{enumerate}
\item $\deg(L)\leq\alpha$, if $\gamma(L)=0$. Moreover, if
  $\deg(L)=\alpha$, there is an $L'\subset V$ such that $V=L\oplus L'$
  and with respect to this decomposition,
$$\gamma=\begin{pmatrix}
          0 & 0\\
	  0 & \gamma'
         \end{pmatrix}$$ with $\gamma'\in H^0(X,L'^{-2}U)$ non-zero;
\item $\deg(L)\leq d/2$, if $\gamma(L)\subset L^\perp U$. Moreover, if $\deg(L)=d/2$, there is $L'\subset V$ such that $V=L\oplus L'$ and with respect to this decomposition,
$$\gamma=\begin{pmatrix}
          0 & \gamma'\\
	  \gamma' & 0
         \end{pmatrix}$$ with $\gamma'\in H^0(X,L^{-1}L'^{-1}U)$ non-zero;
\item $\deg(L)\leq d-\alpha$, if $\gamma(L)\not\subset L^\perp U$. Moreover, if $\deg(L)=d-\alpha$, there is $L'\subset V$ such that $V=L\oplus L'$ and with respect to this decomposition, 
$$\gamma=\begin{pmatrix}
          \gamma' & 0\\
	  0 & 0
         \end{pmatrix}$$ with $\gamma'\in H^0(X,L^{-2}U)$ non-zero.
\end{enumerate}
\end{itemize}
\end{proposition}

\begin{definition}
  Let $(V,\gamma)$ be a quadratic pair. A subbundle $L\subset V$ is \emph{$\alpha$-destabilizing of
    type}:
  \begin{enumerate}
  \item[\textbf{(A)}] if $\deg(L)\geq \alpha$ and $\gamma(L)=0$;
  \item[\textbf{(B)}] if $\deg(L)\geq d/2$ and $\gamma(L)\subset L^{\perp}U$;
\item[\textbf{(C)}] if $\deg(L)\geq d-\alpha$ and $\gamma(L)\not\subset  L^{\perp}U$.
  \end{enumerate}
\end{definition}

\begin{proposition}\label{destab}
 Let $(V,\gamma)$ be a quadratic pair and let $\alpha<d/2$.
\begin{enumerate}
\item There is at most one $\alpha$-destabilizing subbundle $L\subset
  V$ of type \textbf{(A)} and at most one $\alpha$-destabilizing
  subbundle $M\subset V$ of type \textbf{(C)}. Moreover, if such $L$ and $M$ both exist, then $V\cong L\oplus M$.
\item There are at most two distinct $\alpha$-destabilizing subbundles
  $L_1,L_2\subset V$ of type \textbf{(B)}. Moreover, if there exist such distinct $L_1$ and $L_2$, then
  $V\cong L_1\oplus L_2$ and $\gamma(L_1)\subset L_2^{-1} U$.
\item There cannot exist simultaneously $\alpha$-destabilizing
  subbundles of type \textbf{(A)} and \textbf{(B)}.
\item There cannot exist simultaneously $\alpha$-destabilizing
  subbundles of type \textbf{(C)} and \textbf{(B)}.
\end{enumerate}
\end{proposition}
\proof Since $\rk(V) = 2$ and $\gamma$ is holomorphic and non-zero,
there is at most one subbundle $L\subset V$ with $\gamma(L)=0$. This
proves the first statement in $(1)$.  For the second statement in
$(1)$, note that $\deg(M)\geq d-\alpha>d/2=\mu(V)$. Thus the claim
about the destabilizing bundle $M$ follows from the uniqueness of
destabilizing subbundles of ordinary rank $2$ vector bundles
(cf. Proposition $10.38$ of \cite{mukai:2003}). If such $L$ and $M$
both exist, then clearly $L\ncong M$ so the composite $M\to
V\to\Lambda^2VL^{-1}$ is non-zero, and $-\deg(M)+d-\deg(L)\geq 0$. But
$\deg(M)\geq d-\alpha$ and $\deg(L)\geq\alpha$, therefore
$\Lambda^2VL^{-1}\cong M$ and $V\cong L\oplus M$.

The proof of $(2)$ is similar. Let $L_1,L_2\subset V$ be two different
destabilizing subbundles of $(V,\gamma)$ of type \textbf{(B)}. Then
$L_2\cong\Lambda^2VL_1^{-1}$ as before, and $V=L_1\oplus L_2$. In this
case, $L_2^{-1}\cong L_1^\perp$ and $\gamma(L_1)\subset L_2^{-1}
U$. It is clear that there cannot exist a third subbundle satisfying
the same conditions.

The proof of (3) is similar to the proof of the first statement of (1)
because $\gamma(L)=0$ is equivalent to $\gamma(V)\subset L^\perp
U$. The proof of (4) is analogous to the proof of the second statement
in (1), observing that there cannot exist simultaneously a
destabilizing subbundle of $V$ and another subbundle with degree
$d/2$. Indeed, if there is an $M \subset V$ with $\deg(M)>d/2$ and if
$F\subset V$ is different from $M$ then there is a non-zero
homomorphism $F\to \Lambda^2VM^{-1}$, so $\deg(F)\leq d-\deg(M)<d/2$.
\endproof


Recall that $d_U=\deg(U)$ and that $\cN_\alpha(2,d)$ denotes the
moduli space of $\alpha$-polystable $U$-quadratic pairs of rank $2$
and degree $d$.
\begin{proposition}\label{empty above d_U}
\mbox{}
\begin{enumerate}
 \item If $d>d_U$, then $\cN_\alpha(2,d)=\emptyset$ for all $\alpha$.
 \item If $d\leq d_U$, then $\cN_\alpha(2,d)=\emptyset$ for all $\alpha>d/2$.
\end{enumerate}
\end{proposition}
\proof Let $(V,\gamma)$ be a quadratic pair of rank $2$ and degree
$d>d_U$.  If $\rk(\gamma)=2$ (generically), then $\det(\gamma)$ is a
non-zero section of $\Lambda^2V^{-2}U^2$ so $d\leq d_U$. Hence, since $\gamma\neq 0$, we must
have $\rk(\gamma)=1$.  Take any $\alpha$ and suppose moreover that
the pair $(V,\gamma)$ is $\alpha$-semistable. Since $V$ is locally
free, the sheaf $N=\ker(\gamma) \subset V$ is torsion free. For the
same reason, the quotient $V/N\cong \im(\gamma) \subset V^*\otimes U$
is torsion free. Thus $N$ is a line subbundle of $V$. Let
$I\subset V^*$ be such that $IU$ is the saturation of the image sheaf
$\im(\gamma)$.  From the $\alpha$-semistability condition,
\begin{equation}\label{deg kernel}
 \deg(N)\leq\alpha
\end{equation}
and, since $\gamma(I^\perp)=0$,
\begin{equation}\label{deg image}
 \deg(I)\leq\alpha-d.
\end{equation}
On the other hand, $\gamma$ induces a non-zero map of line bundles
$V/N\to IU$, so
\begin{equation}\label{non-zero map}
 -d+\deg(N)+\deg(I)+d_U\geq 0.
\end{equation}
But, from (\ref{deg kernel}) and (\ref{deg image}), we have
\begin{equation}\label{-d+dn+di+du}
 -d+\deg(N)+\deg(I)+d_U<0
\end{equation}
because $d>d_U$ and $\alpha\leq d/2$.  From (\ref{non-zero map}) and
(\ref{-d+dn+di+du}) we conclude that there is no such $(V,\gamma)$ and
this finishes the proof of the first part.

The second part is immediate, since $\alpha\leq
d/2$ is part of the definition of $\alpha$-semistability.
\endproof

This result deals with the cases $d>d_U$ and any $\alpha$, and $d\leq
d_U$ and $\alpha>d/2$. From now on we will restrict ourselves to the
study of $U$-quadratic pairs of type $(2,d)$ with $d<d_U$. When $d=d_U$, the map
$\gamma$ becomes an isomorphism, making this a special case in what
concerns the connected components of the moduli space. In the next remark we give a very brief explanation of this phenomenon, which can be seen as somewhat similar to the difference between
the situations in Propositions \ref{description of N(1,dU/2)} and
\ref{description of Nalpha(1,d) as a cover}.

\begin{remark}
If $d_U$ is odd and $d=d_U$, then it will follow from Proposition \ref{prop:critical values} below (see also Figure 1 in section \ref{sec:critical values}) that $\cN_\alpha(2,d_U)=\emptyset$. 
Assume hence, that $d_U$ is even and that $d=d_U$. In this case $\gamma:V\to V^*\otimes U$ is an isomorphism. If we choose a square root $U'$ of $U$, then $\gamma$ gives rise to a symmetric isomorphism $q:V\otimes U'^*\cong V^*\otimes U'$ i.e. to a non-degenerate quadratic form on the vector bundle $V\otimes U'^*$. Moreover, it can be seen that the $\mathrm{O}(2,\C)$-bundle $(V\otimes U'^*,q)$ is semistable (i.e. the degree of any isotropic subbundle of $V\otimes U'^*$ is less or equal than $0$) if and only if $(V,\gamma)$ is $\alpha$-semistable for any $\alpha\leq d/2$. Hence $\cN_\alpha(2,d_U)$ is isomorphic to the moduli space $\cM_{\mathrm{O}(2,\C)}$ of orthogonal bundles and this gives rise to the existence of extra connected components (cf.\
\cite{bradlow-garcia-prada-gothen:2005,garcia-gothen-mundet:2008
  II,gothen:2001}).  
\end{remark}

\subsection{Deformation theory of quadratic pairs}
\label{sec:deformation-theory}

The deformation theory of a quadratic pair $(V,\gamma)$ is governed by
the following complex of sheaves on $X$ (see, e.g.,
Biswas--Ramanan \cite{biswas-ramanan:1994}):
$$
C^\bullet(V,\gamma):\End(V)\xrightarrow{\rho(\gamma)}S^2V^*\otimes U,
$$
where 
$$\rho(\gamma)(\psi)=-(\psi^t\otimes 1_U)\gamma-\gamma\psi.
$$
In particular, the infinitesimal deformation space of a quadratic pair 
$(V,\gamma)$ is isomorphic to
$\mathbb{H}^1(X,C^\bullet(V,\gamma))$. 
Moreover, one has a long exact
sequence
\begin{equation}
\label{eq:les}
\begin{split}
0&\longrightarrow\mathbb{H}^0(X,C^\bullet(V,\gamma))\longrightarrow
H^0(X,\End(V))\longrightarrow
H^0(X,S^2V^*\otimes U)\longrightarrow\\
&\longrightarrow\mathbb{H}^1(X,C^\bullet(V,\gamma))\longrightarrow H^1(X,\End(V))\longrightarrow
H^1(X,S^2V^*\otimes U)\longrightarrow\\
&\longrightarrow\mathbb{H}^2(X,C^\bullet(V,\gamma))\longrightarrow 0
\end{split}
\end{equation}
where the maps $H^i(X,\End(V))\to H^i(X,S^2V^*\otimes U)$ are induced
by $\rho(\gamma)$.  It is immediate from this long exact sequence that
the infinitesimal automorphism space (defined for general pairs in
\cite{garcia-gothen-mundet:2008}) of a quadratic
pair $(V,\gamma)$ can be canonically identified with
$\mathbb{H}^0(X,C^\bullet(V,\gamma))$.

\begin{definition}\label{def of simple in quadratic pairs}
  A quadratic pair $(V,\gamma)$ is \emph{infinitesimally simple} if
  the vanishing $\mathbb{H}^0(X,C^\bullet(V,\gamma))=0$ holds.
  A quadratic pair $(V,\gamma)$ is \emph{simple} if the group
  $\Aut(V,\gamma)$ of automorphisms of $(V,\gamma)$ is equal to $\{\pm
  1_V\}$.
\end{definition}

The following is a standard fact.
\begin{proposition}
  \label{prop:simple-H2-smooth}\mbox{}
  \begin{enumerate}
  \item An $\alpha$-stable quadratic pair is infinitesimally simple.
  \item An $\alpha$-stable quadratic pair $(V,\gamma)$ represents a
    smooth point in the moduli space if it is simple and
    $\mathbb{H}^2(X,C^\bullet(V,\gamma))=0$.
  \end{enumerate}
\end{proposition}

This motivates the following definition.

\begin{definition}\label{expected dimension}
  The \emph{expected dimension} of $\cN_\alpha(2,d)$ is
  $\dim\mathbb{H}^1(X,C^\bullet(V,\gamma)).$
\end{definition}
Using (\ref{eq:les}), the expected dimension can be calculated as
follows:
\begin{equation}\label{dimension of moduli of pairs}
  \dim\mathbb{H}^1(X,C^\bullet(V,\gamma))
  =\chi(S^2V^*\otimes U)-\chi(\End(V))=3(d_U-d)+g-1.
\end{equation}

\begin{remark}
  If a (local) universal family exists over a component of the moduli
  space, then this component has the expected dimension.  Notice,
  however, that the actual dimension of the moduli space can be strictly smaller than the expected dimension (see
  \cite{bradlow-garcia-prada-gothen:2003} for an example of this
  phenomenon, in the Higgs bundle context).
\end{remark}

\section{Variation of the moduli with the parameter}
\label{sec:variation-of-moduli}

The purpose of this section is to study the variation of the moduli
spaces $\mathcal{N}_{\alpha}(2,d)$ with the stability parameter
$\alpha$.  As in the case of holomorphic triples
\cite{bradlow-garcia-prada-gothen:2004 triples,thaddeus:1994} we have critical
values $\alpha_k$ --- for which the moduli spaces
$\mathcal{N}_{\alpha}(2,d)$ change --- and corresponding \emph{flip
  loci} $\cS_{\alpha_k^\pm}(2,d) \subset
\mathcal{N}_{\alpha^\pm}(2,d)$, where the change takes place. We shall
see that, in contrast with the case of holomorphic triples, there is
no symmetry between $\cS_{\alpha_k^+}(2,d)$ and
$\cS_{\alpha_k^-}(2,d)$. This is due to the non-linear nature of
quadratic pairs.

\subsection{Critical values}\label{sec:critical values}
We begin by determining precisely the critical values of the parameter
$\alpha$. 


\begin{proposition}\label{injectivity}
  If $(V,\gamma)$ is an $\alpha$-semistable pair with
  $\alpha<d-d_U/2$, then generically $\rk(\gamma)=2$.
\end{proposition}
\proof 

Recall that we always have $\gamma\neq 0$. If $\rk(\gamma)=1$,
considering again the line bundles $N=\ker(\gamma)\subset V$ and
$I\subset V^*$ such that $IU$ is the saturation of the image sheaf
$\im(\gamma)$, we have, as in the proof of Proposition \ref{empty above d_U}, that
$$
0\leq -d+\deg(N)+\deg(I)+d_U\leq 2\alpha-2d+d_U,
$$ 
i.e., $\alpha\geq d-d_U/2$.
\endproof

The next result shows that the injectivity parameter $d-d_U/2$ of
Proposition \ref{injectivity} is also a stabilization parameter, in
the sense that after it the moduli spaces $\cN_\alpha(2,d)$, for
different values of $\alpha$, are all isomorphic.

\begin{proposition}\label{stabilization}
If $\alpha_2\leq\alpha_1<d-d_U/2$, then a quadratic pair
  $(V,\alpha)$ is $\alpha_1$-semistable if and only if it is
  $\alpha_2$-semistable, and hence
  $\cN_{\alpha_1}(2,d)\simeq\cN_{\alpha_2}(2,d)$.
\end{proposition}
\proof 

Let $(V,\gamma)\in\cN_{\alpha_1}(2,d)$. Since $\alpha_2\leq\alpha_1$, the existence of an $\alpha_2$-destabilizing subbundle implies that it must be of type \textbf{(A)}, which in turn implies that $\rk(\gamma)=1$ generically. But this is impossible due to Proposition \ref{injectivity}, since $\alpha_1<d-d_U/2$. Hence $\cN_{\alpha_1}(2,d)\subseteq\cN_{\alpha_2}(2,d)$.

Conversely, if $(V,\gamma)\in\cN_{\alpha_2}(2,d)$, then $(V,\gamma)\in\cN_{\alpha_1}(2,d)$ unless there is an $\alpha_1$-destabilizing subbundle of $(V,\gamma)$. Hence $L$ is such that $d-\alpha_1<\deg(L)\leq d-\alpha_2$, and $\gamma(L)\not\subset L^\perp U$, therefore the composite 
$L\to V\xrightarrow{\gamma} V^*\otimes U\to L^{-1}U$ is non-zero. Thus
$$-2\deg(L)+d_U\geq 0.$$ On the other hand, $d-\alpha_1<\deg(L)$ together with $\alpha_1<d-d_U/2$, gives
$$-2\deg(L)+d_U<0.$$
It follows that no such line subbundle $L$ can exist.
\endproof

{}From the definition of $\alpha$-semistability and from the previous
proposition, the following is immediate.

\begin{proposition}\label{prop:critical values}
The critical values of $U$-quadratic pairs of type $(2,d)$ are the elements of the following set:
$$\left\{d/2\right\}\cup\left\{[d/2]+k\suchthat\,k\in\left\{d-[d/2]-[d_U/2],\ldots,0\right\}\right\}.$$
Moreover, on each open interval between consecutive critical values,
$$\left(\,[d/2]+k,\min\left\{d/2,[d/2]+k+1\right\}\,\right)$$
the $\alpha$-semistability condition is the same; hence the corresponding moduli spaces are isomorphic.
\end{proposition}

\begin{notation}\label{formula alphak}
  For each integer $d-[d/2]-[d_U/2]\leq k\leq 0$,
  we define
  $$\alpha_k=[d/2]+k.$$
  Also, let 
  $$
  \alpha_M=d/2\quad\text{and}\quad
  \alpha_m=\alpha_{d-[d/2]-[d_U/2]}=d-[d_U/2],
  $$
  and let $\alpha_k^+$ denote the value of any parameter between the
  critical values $\alpha_k$ and $\alpha_{k+1}$, and let $\alpha_k^-$
  denote the value of any parameter between the critical values
  $\alpha_{k-1}$ and $\alpha_k$.
\end{notation}

Proposition~\ref{stabilization} means that we can write without ambiguity 
$$\cN_{\alpha_k^+}(2,d)$$
for the moduli space of $\alpha_k^+$-semistable $U$-quadratic pairs of
rank $2$ and degree $d$, for any $\alpha$ between the critical values
$\left[d/2\right]+k$ and
$\min\left\{d/2,\left[d/2\right]+k+1\right\}$. Note that, with this
notation, we always have $\cN_{\alpha_k^+}(2,d)=\cN_{\alpha_{k+1}^-}(2,d)$.

\begin{figure}[htb]
\centering
\includegraphics{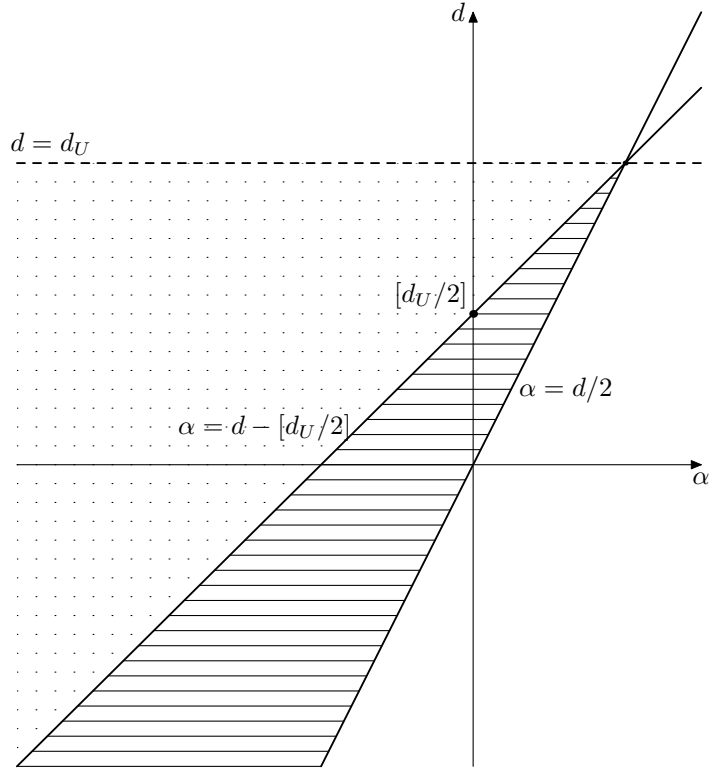}
\caption{Variation regions of $\cN_\alpha(2,d)$. Above the line
  $d=d_U$, $\cN_\alpha(2,d)=\emptyset$ as well as on the right of the
  line $\alpha=d/2$. Also, $\cN_\alpha(2,d_U)=\emptyset$ whenever $d_U$ is odd. The dotted region, on the left of the line
  $\alpha=d-[d_U/2]$, is the region where there are no critical
  values, hence there is no change of the moduli spaces and, also,
  $\gamma$ is non-degenerate. The critical values lie in the region
  between the lines $\alpha=d-[d_U/2]$ and $\alpha=d/2$.}
  \label{fig:alpha-d}
\end{figure}

The information obtained so far on the variation of $\mathcal{N}_{\alpha}(2,d)$
with $\alpha$ and $d$ is summarized in
Figure~\ref{fig:alpha-d}. 

\subsection{Flip loci}\label{Flip loci}

We shall now study what are the differences between moduli spaces of
$U$-quadratic pairs of type $(2,d)$, which are separated by a critical
value of the parameter $\alpha$.

\begin{definition}
  For each $k\in\left\{d-[d/2]-[d_U/2],\ldots,0\right\}$, let
  $\cS_{\alpha_k^+}(2,d)$ be the set of pairs of degree $d$ which are
  $\alpha_k^+$-semistable but $\alpha_k^-$-unstable, i.e.,
  $$
  \cS_{\alpha_k^+}(2,d)=\left\{(V,\gamma)\in\cN_{\alpha_k^+}(2,d)
    \suchthat(V,\gamma)\notin\cN_{\alpha_k^-}(2,d)\right\}.
  $$
  Similarly, define $\cS_{\alpha_k^-}(2,d)$ to be the set of pairs of
  degree $d$ which are $\alpha_k^-$-semistable but
  $\alpha_k^+$-unstable,
  $$
  \cS_{\alpha_k^-}(2,d)=\left\{(V,\gamma)\in\cN_{\alpha_k^-}(2,d)
    \suchthat(V,\gamma)\notin\cN_{\alpha_k^+}(2,d)\right\}.
  $$
  The spaces $\cS_{\alpha_k^\pm}(2,d)$ are called the \emph{flip loci}
  for the critical value $\alpha_k$.
\end{definition}

As a direct consequence of this definition, we have
\begin{equation}\label{N--S-=N+-S+}
  \cN_{\alpha_k^+}(2,d)\setminus\cS_{\alpha_k^+}(2,d)
  =\cN_{\alpha_k^-}(2,d)\setminus\cS_{\alpha_k^-}(2,d).
\end{equation}

\begin{proposition}
  Any quadratic pair $(V,\gamma)$ in $\cS_{\alpha_k^{\pm}}(2,d)$ is
  $\alpha_k^{\pm}$-stable.  Hence, for $\alpha_k\neq\alpha_M$,
  \begin{equation}\label{S+ inside stable locus}
    \cS_{\alpha_k^+}(2,d)=\left\{(V,\gamma)\in\cN^s_{\alpha_k^+}(2,d)
      \suchthat(V,\gamma)\notin\cN_{\alpha_k^-}(2,d)\right\}
  \end{equation}
  and
  \begin{equation}\label{S- inside stable locus}
    \cS_{\alpha_k^-}(2,d)=\left\{(V,\gamma)\in\cN^s_{\alpha_k^-}(2,d)
      \suchthat(V,\gamma)\notin\cN_{\alpha_k^+}(2,d)\right\}
  \end{equation}
  where $\cN^s_{\alpha_k^\pm}(2,d)$ stands for the stable locus of
  $\cN_{\alpha_k^\pm}(2,d)$.
\end{proposition}

\begin{proof}
  If $(V,\gamma)$ is strictly $\alpha_k^\pm$-semistable then, since
  $\alpha_k^\pm$ is a generic value, the destabilizing subbundle must
  be of type \textbf{(B)}. Since such a subbundle is destabilizing for
  all values of $\alpha$, (3) and (4) of Proposition~\ref{destab}
  imply that there are no destabilizing subbundles of type
  \textbf{(A)} or \textbf{(C)}. The conclusion is now immediate from
  the definition of $\cS_{\alpha_k^\pm}(2,d)$.
\end{proof}

\begin{proposition}\label{S+- is a subvariety}
 Let $\alpha_k\neq\alpha_M$. Then:
\begin{enumerate}
 \item $\cS_{\alpha_k^+}(2,d)$ is a subvariety of $\cN^s_{\alpha_k^+}(2,d)$.
 \item $\cS_{\alpha_k^-}(2,d)$ is a subvariety of $\cN^s_{\alpha_k^-}(2,d)$.
\end{enumerate}
\end{proposition}
\proof From (\ref{S+ inside stable locus}),
$\cS_{\alpha_k^+}(2,d)\subset\cN^s_{\alpha_k^+}(2,d)$. From
\cite{gomez-sols:2000} we know that there is a (universal) family of
quadratic pairs parametrized by $\cN^s_{\alpha_k^+}(2,d)$. By
definition, the restriction of this family to $\cS_{\alpha_k^+}(2,d)$
parametrizes the pairs which are not $\alpha_k^-$-semistable. Since
$\alpha_k^-$-semistability is an open condition (cf.\ Proposition~3.1
of \cite{gomez-sols:2000}), it follows that
$\cS_{\alpha_k^+}(2,d)$ is Zariski closed in
$\cN^s_{\alpha_k^+}(2,d)$. This proves (1). The proof of (2) is the
same, but now using (\ref{S- inside stable locus}).
\endproof

\begin{remark}\label{S+- is a subvariety of N}
  In the next two sections we shall see that the flip loci
  $\cS_{\alpha_k^\pm}(2,d)$ are compact and therefore also closed in
  $\cN_{\alpha_k^\pm}(2,d)$. Hence the $\cS_{\alpha_k^\pm}(2,d)$ are in
  fact subvarieties of $\cN_{\alpha_k^\pm}(2,d)$.
\end{remark}

\subsection{The flip locus \texorpdfstring{$\cS_{\alpha_k^+}(2,d)$}{Salpha+(2,d)}}

From (2) of Proposition \ref{empty above d_U},
$\cS_{\alpha_M^+}(2,d)=\emptyset$, so we shall study the flip loci
$\cS_{\alpha_k^+}(2,d)$ for the other critical values.

\begin{proposition}\label{ext S+}
  Let $(V,\gamma)\in\cN_{\alpha_k^+}(2,d)$ with $\alpha_k<d/2$. If
  $(V,\gamma)\in\cS_{\alpha_k^+}(2,d)$, then $V$ is a non-trivial
  extension
  $$0\longrightarrow L\longrightarrow V\longrightarrow M\longrightarrow 0$$
  where $L\subset V$ is a line bundle such that $\deg(L)=\alpha_k$ and
  $\gamma(L)=0$. Moreover, $\gamma$ induces $\gamma'\in
  H^0(X,M^{-2}U)$ such that the quadratic pair $(M,\gamma')$ of
  type $(1,d-\alpha_k)$ is $(d-\alpha_k^+)$-stable.
\end{proposition}

\proof 
Let $(V,\gamma)\in\cS_{\alpha_k^+}(2,d)$. Then it must be strictly
$\alpha_k$-semistable and, from the definition of
$\cS_{\alpha_k^+}(2,d)$ and Proposition \ref{simplified notion of polystability}, the destabilizing subbundle must be an $L\subset V$
such that $\gamma(L)=0$ and $\deg(L)=\alpha_k$.

Write
\begin{equation}\label{ext+}
0\longrightarrow L\longrightarrow V\longrightarrow\Lambda^2VL^{-1}\longrightarrow 0
\end{equation}
and define $M=\Lambda^2VL^{-1}$.  If we had $V=L\oplus M$, then $M$
would be an $\alpha_k^+$-destabilizing subbundle of $(V,\gamma)$ (of
type \textbf{(C)}), which is not possible. The extension (\ref{ext+})
is thus non-trivial.

Using the symmetry of $\gamma$ and the fact that $\gamma(L)=0$, we see
that $\gamma$ induces a map $\gamma':M\to M^{-1}U$ and hence we obtain
the pair $(M,\gamma')$ of type $(1,d-\alpha_k)$. From Definition
\ref{sspairs1}, it is clearly $(d-\alpha_k^+)$-stable.
\endproof


\begin{proposition}\label{S+}
Let $\alpha_k<d/2$. There is a morphism $$\cS_{\alpha_k^+}(2,d)\longrightarrow\cN_{d-\alpha_k^+}(1,d-\alpha_k)\times\Jac^{\alpha_k}(X)$$ with fibre isomorphic to $\mathbb{P}^{d-2\alpha_k+g-2}$.
\end{proposition}
\proof
From Proposition \ref{ext S+}, we see that there is a map
$$p:\cS_{\alpha_k^+}(2,d)\longrightarrow\cN_{d-\alpha_k^+}(1,d-\alpha_k)\times\Jac^{\alpha_k}(X)$$ defined by $$p(V,\gamma)=((M,\gamma'),L).$$
where $L$ is the destabilizing subbundle and $M$ is the quotient bundle, $M\cong\Lambda^2VL^{-1}$.

Let $$((M,\gamma'),L)\in\cN_{d-\alpha_k^+}(1,d-\alpha_k)\times\Jac^{\alpha_k}(X).$$
The fibre of $p$ over
$((M,\gamma'),L)$ is given by the isomorphism classes of non-trivial extensions of $M$ by $L$.
Indeed, if $V$ is such an extension then, defining $$\gamma=(\pi^t\otimes 1_U)\gamma'\pi$$ where $\pi:V\to M$ is the projection, we obtain a quadratic pair $(V,\gamma)$.
This pair is strictly $\alpha_k$-semistable and $\alpha_k^-$-unstable (with $L$ being the destabilizing subbundle) and as we go from $\cN_{\alpha_k^-}(2,d)$ to $\cN_{\alpha_k^+}(2,d)$, then $(V,\gamma)$ 
gets $\alpha_k^+$-stable unless $V$ has also a destabilizing subbundle $M'$ such that
$\deg(M')=d-\alpha_k$. But then $M'\cong M$ and $V=L\oplus M$, which contradicts the non-triviality of the extension $V$.

The fibre of $p$ over $((M,\gamma'),L)$ is then the space
 $\mathbb{P}\Ext^1(M,L)\cong\mathbb{P}H^1(X,M^{-1}L)$. Since
$\alpha_k<d/2$, $\deg(M^{-1}L)=2\alpha_k-d<0$, so $H^0(X,M^{-1}L)=0$ and 
$$\dim H^1(X,M^{-1}L)=d-2\alpha_k+g-1>0.$$
Hence $p$ is surjective, with fibre isomorphic to $\mathbb{P}^{d-2\alpha_k+g-2}$.

It remains to check that $p$ is a morphism. For that we proceed as follows.

Let $p_X:\cN_{d-\alpha_k^+}(1,d-\alpha_k)\times X\to X$ be the projection. From Remark \ref{degenerated cases} and Theorem I of \cite{gomez-sols:2000}, there is a universal $p_X^*U$-quadratic pair $(\mathcal L_1,\pmb{\gamma}')$ over $\cN_{d-\alpha_k^+}(1,d-\alpha_k)\times X$. On the other hand, we have the Poincaré line bundle $\mathcal L_2$ over $\Jac^{\alpha_k}(X)\times X$.
Let $$\mathrm{pr}_{13}:\cN_{d-\alpha_k^+}(1,d-\alpha_k)\times\Jac^{\alpha_k}(X)\times X\longrightarrow\cN_{d-\alpha_k^+}(1,d-\alpha_k)\times X$$
$$\mathrm{pr}_{23}:\cN_{d-\alpha_k^+}(1,d-\alpha_k)\times\Jac^{\alpha_k}(X)\times X\longrightarrow\Jac^{\alpha_k}(X)\times X$$
and $$\mathrm{pr}_{12}:\cN_{d-\alpha_k^+}(1,d-\alpha_k)\times\Jac^{\alpha_k}(X)\times X\longrightarrow\cN_{d-\alpha_k^+}(1,d-\alpha_k)\times\Jac^{\alpha_k}(X)$$ be the projections.
Consider the first direct image sheaf $R^1\mathrm{pr}_{12*}(\mathrm{pr}_{13}^*\mathcal{L}_1^{-1}\otimes\mathrm{pr}_{23}^*\mathcal{L}_2)$. This sheaf is locally free since its fibres have constant dimension (because $\alpha_k<d/2$).

If $$\widetilde{\mathcal S}((\mathcal L_1,\pmb{\gamma}'),\mathcal L_2)=\mathbb{P}R^1\mathrm{pr}_{12*}(\mathrm{pr}_{13}^*\mathcal{L}_1^{-1}\otimes\mathrm{pr}_{23}^*\mathcal{L}_2)$$ there is then a morphism $\widetilde{\mathcal S}((\mathcal L_1,\pmb{\gamma}'),\mathcal L_2)\to\cN_{d-\alpha_k^+}(1,d-\alpha_k)\times\Jac^{\alpha_k}(X)$. Moreover, in a similar manner to \cite{lange:1983} (see also Proposition $3.2$ of \cite{thaddeus:1994} and Proposition $5.10$ of \cite{garcia-prada-gothen-munoz:2007}), one sees that $\widetilde{\mathcal S}((\mathcal L_1,\pmb{\gamma}'),\mathcal L_2)$ is base of a family parametrizing all $\alpha_k^+$-semistable $U$-quadratic pairs over $X$ which are $\alpha_k^+$-stable but $\alpha_k^-$-unstable. Hence, from the universal property of the coarse moduli space $\cN_{\alpha_k^+}(2,d)$, there is a morphism $\widetilde{\mathcal S}((\mathcal L_1,\pmb{\gamma}'),\mathcal L_2)\to\cN_{\alpha_k^+}(2,d)$ which factors through $\mathcal S_{\alpha_k^+}(2,d)$ and yields an isomorphism $\widetilde{\mathcal S}((\mathcal L_1,\pmb{\gamma}'),\mathcal L_2)\cong\mathcal S_{\alpha_k^+}(2,d)$ such that the following diagram commutes:
$$\xymatrix{\widetilde{\mathcal S}((\mathcal L_1,\pmb{\gamma}'),\mathcal L_2)\ar[r]^{\cong}\ar[rd]&\mathcal S_{\alpha_k^+}(2,d)\ar[d]^{p}\\
&\cN_{d-\alpha_k^+}(1,d-\alpha_k)\times\Jac^{\alpha_k}(X).}$$ So $p$ is a morphism and the result follows.
\endproof

Since $\cN_{d-\alpha_k^+}(1,d-\alpha_k)$ is compact, one concludes from this proposition that $\cS_{\alpha_k^+}(2,d)$ is compact as well. It follows from Proposition \ref{S+- is a subvariety} (see also Remark \ref{S+- is a subvariety of N}) that it is a subvariety of $\cN_{\alpha_k^+}(2,d)$.

From the previous proposition and from Proposition \ref{description of Nalpha(1,d) as a cover} (in the case of $\cS_{\alpha_m^+}(2,d)$ use instead Proposition \ref{description of N(1,dU/2)}), we have:
\begin{corollary}\label{Salphak+(2,d)}
For every $\alpha_k<d/2$, $\dim\cS_{\alpha_k^+}(2,d)=d_U-d+2g-2$.
\end{corollary}

\subsection{The flip locus \texorpdfstring{$\cS_{\alpha_k^-}(2,d)$}{Salpha-(2,d)}}\label{subsection flip S-}

Now we turn our attention to the other flip loci, $\cS_{\alpha_k^-}(2,d)$. As in the case of $\cS_{\alpha_k^+}(2,d)$, the behaviour of $\cS_{\alpha_k^-}(2,d)$ depends on whether $\alpha_k=\alpha_M=d/2$ or not. On the other hand, for $\alpha_k<d/2$, the description of $\cS_{\alpha_k^-}(2,d)$ is more involved than that of $\cS_{\alpha_k^+}(2,d)$, with several difficulties appearing due to the fact that the $\alpha_k$-destabilizing subbundle is of type \textbf{(C)}.

Let us begin by studying $\cS_{\alpha_M^-}(2,d)$ and see why it is a separate case.
Indeed, $$\cS_{\alpha_M^-}(2,d)=\cN_{\alpha_M^-}(2,d)$$ and one cannot compare this flip locus with the others because in this extreme case the cause of the destabilization after $d/2$ is not related with subbundles of $V$.

We have that $\alpha_M^-$ is any value in the open interval $([d/2],d/2)$ if $d$ is odd, or $(d/2-1,d/2)$ if $d$ is even. So we can write $\alpha_M^-=d/2-\epsilon$ for sufficiently small $\epsilon>0$.

\begin{lemma}
 If a $U$-quadratic pair $(V,\gamma)$ of type $(2,d)$ is $\alpha_M^-$-semistable then $V$ is semistable.
\end{lemma}
\proof
Suppose $(V,\gamma)$ is $\alpha_M^-$-semistable and let $L\subset V$ be a line subbundle. Then:
\begin{itemize}
 \item $\deg(L)<\alpha_M^-\leq d/2$ if $\gamma(L)=0$ (note that $\alpha_M^-$ is not a critical value, so we could not have $\deg(L)=\alpha_M^-$);
 \item $\deg(L)\leq d/2$ if $\gamma(L)\subset L^\perp U$;
 \item $\deg(L)< d-\alpha_M^-=d/2+\epsilon$ if $\gamma(L)\not\subset L^\perp U$ (again we could not have $\deg(L)=d-\alpha_M^-$). Since we can take $\epsilon>0$ as small as wanted, it follows that $\deg(L)\leq d/2$.
\end{itemize}
In any case, $\deg(L)\leq d/2$, and thus $V$ is semistable.
\endproof


\begin{proposition}\label{fibre d/2}
Let $\mathcal{M}(2,d)$ be the moduli space of rank $2$ semistable vector bundles over $X$. There is a map $\pi:\cN_{\alpha_M^-}(2,d)\to\mathcal{M}(2,d)$ which, if $d_U-d>g-1$, is surjective and the fibre over a stable vector bundle $V$ is $\mathbb{P}H^0(X,S^2V^*\otimes U)$.
\end{proposition}
\proof
Using the previous lemma, define $\pi:\cN_{\alpha_M^-}(2,d)\to\mathcal{M}(2,d)$ as $\pi(V,\gamma)=V$.

For the map to be surjective, given any semistable vector bundle $V$ in $\cM(2,d)$, there must exist a non-zero holomorphic section $\gamma$ of $S^2V^*\otimes U$ such that $(V,\gamma)$ is $\alpha_M^-$-semistable.
Suppose that $V$ is stable. Since $d_U-d>g-1$, we have $\chi(S^2V^*\otimes U)>0$, where $\chi$ is the Euler characteristic, hence $H^0(X,S^2V^*\otimes U)\neq 0$ and this yields a pair $(V,\gamma)$.
For any $L\subset V$, $\deg(L)<d/2$. Hence $\deg(L)\leq d/2-1$ if $d$ is even or $\deg(L)\leq [d/2]$.
In any case, $\deg(L)\leq\alpha_M^-$, hence $(V,\gamma)$ is $\alpha_M^-$-stable. The image of $\pi$ contains therefore the open dense subspace of stable vector bundles. Since $\cN_{\alpha_M^-}(2,d)$ is compact (cf. Theorem I of \cite{gomez-sols:2000}) it follows that $\pi$ is surjective.

To compute the fibre over a stable vector bundle, we only have to note that $(V,\gamma)\cong(V,\gamma')$ if and only if $\gamma=\theta\gamma'$ for some $\theta\in\C^*$. If $(V,\gamma)\cong(V,\gamma')$, then there is an automorphism $\lambda:V\to V$ such that $\gamma'\lambda=((\lambda^t)^{-1}\otimes 1_U)\gamma$. But, as $V$ is stable, $\lambda$ is a non-zero scalar so $\gamma=\lambda^2\gamma'$.
On the other hand, if $\gamma=\theta\gamma'$ for some $\theta\in\C^*$ then the scalar automorphism of $V$ given by $\sqrt{\theta}\in\C^*$ is an isomorphism between $(V,\gamma)$ and $(V,\gamma')$.
\endproof


We now move on to the description of the flip loci $\cS_{\alpha_k^-}(2,d)$ with $\alpha_k<d/2$.

If $(V,\gamma)\in\cS_{\alpha_k^-}(2,d)$, then it is $\alpha_k^-$-stable and $\alpha_k^+$-unstable hence strictly $\alpha_k$-semistable. The destabilizing subbundle must be a line subbundle $M\subset V$ such that $$\deg(M)=d-\alpha_k$$ and $$\gamma(M)\not\subset M^\perp U.$$
Therefore $\gamma$ induces a non-zero holomorphic map
\begin{equation}\label{gamma'}
\gamma':M\longrightarrow V\stackrel{\gamma}{\longrightarrow} V^*\otimes U\longrightarrow M^{-1}U
\end{equation}
i.e. $$\gamma'\in H^0(X,M^{-2}U).$$
The description of the flip loci $\cS_{\alpha_k^-}(2,d)$, with $\alpha_k<d/2$, will be done by carrying out a detailed analysis of this information.

Write $V$ as an extension
\begin{equation}\label{ext-}
0\longrightarrow M\longrightarrow V\longrightarrow L\longrightarrow 0
\end{equation}
where $L=\Lambda^2VM^{-1}$.

What we have is already enough to describe $\cS_{\alpha_m^-}(2,d)$ for $U$-quadratic pairs such that $d_U$ is even.

\begin{proposition}\label{description of pairs in S-alpham}
Suppose that $d_U$ is even. If $(V,\gamma)\in\cS_{\alpha_m^-}(2,d)$ then the extension (\ref{ext-}) is trivial, $$V=M\oplus L,$$ 
and with respect to this decomposition,
$$\gamma=\begin{pmatrix}
     \gamma' & 0 \\
     0 & \gamma''
\end{pmatrix}$$
where $\gamma'$ is defined in (\ref{gamma'}) and $\gamma''\in H^0(X,L^{-2}U)\setminus\{0\}$.
\end{proposition}
\proof
If $(V,\gamma)\in\cS_{\alpha_m^-}(2,d)$, then $\deg(M)=d-\alpha_m=d_U/2$, thus $\deg(M^{-1}U)=d_U/2$.

Consider the map $$\varphi:V\longrightarrow M^{-1}U$$ defined by $$\varphi=(i^t\otimes 1_U)\gamma,$$ where $i:M\hookrightarrow V$ is the inclusion. Since $\gamma(V)\not\subset M^{\perp}U$, then $\rk(\varphi)=1$ generically.

Consider the line subbundle $N=\ker(\varphi)$ of $V$.
We have the induced non-zero map $$\widetilde\varphi:V/N\longrightarrow M^{-1}U$$ i.e., $\widetilde\varphi\in H^0(X,(V/N)^{-1}M^{-1}U)$, hence 
$$\deg(N)\geq d-\deg(M^{-1}U)=d-d_U/2=\deg(L).$$
On the other hand, since $N\ncong M$, we have a non-zero map $N\to L$, so 
$$\deg(N)\leq\deg(L).$$

We conclude that $\deg(N)=\deg(L)$ and that the map $N\to L$ is an isomorphism, $L\cong N$, from which follows that extension (\ref{ext-}) is trivial: $$V=M\oplus L.$$
Since $L\cong N=\ker(\varphi)$, we have that $\gamma(L)\subset M^\perp U\cong L^{-1}U$, thus the form of $\gamma$ with respect to the decomposition $V=M\oplus L$ is
$$\gamma=\begin{pmatrix}
     \gamma' & 0 \\
     0 & \gamma''
\end{pmatrix}.$$ $\gamma''\neq 0$ otherwise $L$ would be an $\alpha_m^-$-destabilizing subbundle of $(V,\gamma)$ of type \textbf{(A)}, contradicting the assumption $(V,\gamma)\in\cS_{\alpha_m^-}(2,d)$ (or, alternatively, because of Proposition \ref{injectivity}).
\endproof

\begin{corollary}\label{S-alpham}
If $U$ has even degree $d_U$, $\cS_{\alpha_m^-}(2,d)$ is isomorphic to $$\cN_{d_U/2^-}(1,d_U/2)\times\cN_{(d-d_U/2)^-}(1,d-d_U/2).$$
\end{corollary}
\proof
Given $(V,\gamma)\in\cS_{\alpha_m^-}(2,d)$, the pair $(V,\gamma)$ determines and is determined by the pairs $(M,\gamma')$ and $(L,\gamma'')$ obtained in the previous proposition. These are $(d-\alpha_m^+)$-stable and $\alpha_m^-$-stable, respectively,
therefore the map $(V,\gamma)\mapsto((M,\gamma'),(L,\gamma''))$ is an isomorphism between $\cS_{\alpha_m^-}(2,d)$ and $\cN_{d-\alpha_m^+}(1,d-\alpha_m)\times\cN_{\alpha_m^-}(1,\alpha_m)$.
\endproof

So, from Proposition \ref{description of N(1,dU/2)}, in this case $\cS_{\alpha_m^-}(2,d)$ is isomorphic to $$S\times\cN_{(d-d_U/2)^-}(1,d-d_U/2)$$ where $S$ is the set of square roots of $U$.

\vspace{0,5cm}
Now we pass to the analysis of $\cS_{\alpha_k^-}(2,d)$, with $\alpha_k\neq\alpha_m,\alpha_M$ if $d_U$ is even or just $\alpha_k\neq\alpha_M$ if $d_U$ is odd. We start by noticing some constrains of its elements.

In the cases we are now considering, the map $\gamma':M\to M^{-1}U$ as defined in (\ref{gamma'}) is not an isomorphism. Let $D$ be its divisor
\begin{equation}\label{defD}
D=\divisor(\gamma')
\end{equation}
and consider the structure sheaf $\mathcal{O}_D$ of $X$ restricted to $D$ (or structure sheaf of the scheme $D$).

\begin{remark}\label{restriction to a divisor}
If $D=\sum_{i=1}^m n_ip_i$, then, choosing a local coordinate $z_i$ centred at $p_i$, a global section of $\mathcal{O}_D$ can be written as $\sum_{i=1}^m f_i(z)$ where $f_i(z)=\sum_{k=0}^{n_i-1}a_kz_i^k$. One has then a short exact sequence of sheaves
\begin{equation}\label{ses str ring supp on D}
0\longrightarrow\mathcal{O}(-D)\longrightarrow\mathcal{O}\stackrel{r(D)}{\longrightarrow}\mathcal{O}_D\longrightarrow 0
\end{equation} where, for each open $U$ such that there is only one point $p\in\cSupp(D)$ in $U$,
\begin{equation}\label{def of r(D)}
r(D)(U)(s)=r(D)(U)\left(\sum_{k=0}^\infty a_kz^k\right)=\sum_{k=0}^{D(p)-1}a_kz^k,
\end{equation}
for $s\in\mathcal{O}(U)$ such that, in a local coordinate $z$ centred at $p$, $s(z)=\sum_{k=0}^\infty a_kz^k$.
\end{remark}

\begin{proposition}\label{prop:theta-gamma}
  There is a well defined section
  \begin{equation}\label{definition of q}
\theta_\gamma\in H^0(D,M^{-1}L^{-1}U).
  \end{equation}
  given by restriction of $\gamma |_M$ to $D$.
\end{proposition}
\proof
For any sheaf $\mathcal{F}$, write $\mathcal{F}|_D$ for $\mathcal{F}\otimes\mathcal{O}_D$. From (\ref{ses str ring supp on D}), we obtain the short exact sequence of sheaves
\begin{equation}
0\longrightarrow V^*\otimes M^{-1}U(-D)\longrightarrow V^*\otimes M^{-1}U\stackrel{r(D)}{\longrightarrow}V^*\otimes M^{-1}U|_D\longrightarrow 0
\end{equation} and we have a map, which we still denote by $r(D)$,
$$r(D):H^0(X,V^*\otimes M^{-1}U)\longrightarrow H^0(D,V^*\otimes M^{-1}U).$$ 
Now, $\gamma|_M\in H^0(X,V^*\otimes M^{-1}U)$ so consider
$r(D)(\gamma|_M)\in H^0(D,V^*\otimes M^{-1}U)$. But, since $D=\divisor(\gamma')$, we have in fact that $r(D)(\gamma|_M)\in H^0(D,M^{-1}L^{-1}U)$, so we define
$$\theta_\gamma=r(D)(\gamma|_M)\in H^0(D,M^{-1}L^{-1}U),$$ as claimed.
\endproof

If we also denote by $r(D)$ the map in $H^0$ of the restriction $\Lambda^2V^{-2}U^2\to\Lambda^2V^{-2}U^2|_D$, we see that
\begin{equation}\label{det=bargamma2}
 r(D)(\det(\gamma))=-\theta_\gamma^2\in H^0(D,M^{-2}L^{-2}U^2).
\end{equation}
This section $\theta_\gamma\in H^0(D,M^{-1}L^{-1}U)$, obtained in the previous proposition, will be very important in the description of $\cS_{\alpha_k^-}(2,d)$ and the next result is a first instance of this.

\begin{proposition}\label{q=0 implies ...}
 Let $(V,\gamma)\in\cS_{\alpha_k^-}(2,d)$. Then $\theta_\gamma=0$ if and only if extension (\ref{ext-}) is trivial, $$V=M\oplus L,$$ 
and with respect to this decomposition, 
$$\gamma=\begin{pmatrix}
     \gamma' & 0 \\
     0 & \gamma''
\end{pmatrix},$$
where $\gamma'$ is defined in (\ref{gamma'}) and $\gamma''\in H^0(X,L^{-2}U)\setminus\{0\}$. 
\end{proposition}
\proof
As in the proof of Proposition \ref{description of pairs in S-alpham}, consider the map $\varphi:V\to M^{-1}U$ given by $\varphi=(i^t\otimes 1_U)\gamma$, and its kernel $N$.

We have the induced non-zero map $$\widetilde\varphi:V/N\longrightarrow M^{-1}U$$ i.e., $\widetilde\varphi\in H^0(X,(V/N)^{-1}M^{-1}U)$ and let $$\widetilde D=\divisor(\widetilde\varphi)$$ so that 
\begin{equation}\label{deg widetilde D}
 \deg(\widetilde D)=\deg(N)-2d+\alpha_k+d_U.
\end{equation}

Let $p$ be any point in $\cSupp(D)$ and choose a local coordinate $z$ of $X$ centred at $p$. Locally, we can write 
$V=M\oplus L$ and, with respect to this decomposition,
$$\gamma(z)=\begin{pmatrix}
     f_1(z) & f_2(z) \\
     f_2(z) & f_3(z)
\end{pmatrix}$$
hence $$\varphi(z)=\Big( f_1(z) \hspace{.4cm} f_2(z)\Big)$$ so
\begin{equation}\label{D min}
	\widetilde{D}(p)\geq\min\{\ord_0f_1(z),\ord_0f_2(z)\}=\min\{D(p),\text{ord}_0f_2(z)\}. 
\end{equation}
Since $\theta_\gamma=0$, then 
\begin{equation}\label{theta(p)=0}
\theta_\gamma(p)=0
\end{equation}
but, by the definition of $\theta_\gamma$ in (\ref{definition of q}), $$\theta_\gamma(p)=r(D)(U)(f_2(z))$$ and from this, (\ref{def of r(D)}) and (\ref{theta(p)=0}), we see that $\text{ord}_0f_2(z)\geq D(p)$. It follows from
(\ref{D min}) that 
\begin{equation}\label{widetilde D geq D}
\widetilde{D}\geq D
\end{equation}
so $$\deg(\widetilde D)\geq\deg(D)=-2\deg(M)+d_U=-2d+2\alpha_k+d_U.$$
From this and (\ref{deg widetilde D}), we see that
\begin{equation}\label{deg N>=deg L}
\deg(N)\geq\alpha_k=\deg(L). 
\end{equation}
On the other hand, since $N\ncong M$, we have a non-zero map $N\to L$ so (\ref{deg N>=deg L}) implies $L\cong N$. Extension (\ref{ext-}) is hence trivial: $$V=M\oplus L.$$
From $L\cong N=\ker(\varphi)$, it follows that $\gamma(L)\subset M^\perp U\cong L^{-1}U$, thus the form of $\gamma$ with respect to the decomposition $V=M\oplus L$ is
$$\gamma=\begin{pmatrix}
     \gamma' & 0 \\
     0 & \gamma''
\end{pmatrix}.$$ $\gamma''\neq 0$ otherwise $L$ would be an $\alpha_k^-$-destabilizing subbundle of $(V,\gamma)$ of type \textbf{(A)}, contradicting the assumption $(V,\gamma)\in\cS_{\alpha_k^-}(2,d)$.
\endproof

We shall write $\cS_{\alpha_k^-}(2,d)$ as a disjoint union
\begin{equation}\label{dijunion}
\cS_{\alpha_k^-}(2,d)=\cS_{\alpha_k^-}^0(2,d)\sqcup\cS_{\alpha_k^-}^1(2,d)
\end{equation}
where
\begin{itemize}
 \item $\cS_{\alpha_k^-}^0(2,d)$ is the space of pairs in $\cS_{\alpha_k^-}(2,d)$ with $\theta_\gamma=0$;
 \item $\cS_{\alpha_k^-}^1(2,d)$ is the space of pairs in $\cS_{\alpha_k^-}(2,d)$ with $\theta_\gamma\neq 0$.
\end{itemize}

Let us now study each of the spaces $\cS_{\alpha_k^-}^0(2,d)$ and $\cS_{\alpha_k^-}^1(2,d)$.

\subsubsection{}
From Proposition \ref{q=0 implies ...}, $\cS_{\alpha_k^-}^0(2,d)$ is precisely the space of pairs in $\cS_{\alpha_k^-}(2,d)$ such that extension (\ref{ext-}) is trivial and $\gamma$ has the given form. So we have the following corollary.

\begin{corollary}\label{S-alphak0}
Let $\alpha_k\neq\alpha_M$. Then $\cS_{\alpha_k^-}^0(2,d)$ is isomorphic to $$\cN_{d-\alpha_k^+}(1,d-\alpha_k)\times\cN_{\alpha_k^-}(1,\alpha_k).$$
\end{corollary}
\proof
Given $(V,\gamma)\in\cS_{\alpha_k^-}^0(2,d)$, the pair $(V,\gamma)$ determines and is determined by the pairs $(M,\gamma')$ and $(L,\gamma'')$ obtained in the previous proposition. These are $(d-\alpha_k^+)$-stable and $\alpha_k^-$-stable, respectively, therefore the map $(V,\gamma)\mapsto((M,\gamma'),(L,\gamma''))$ is an isomorphism between $\cS_{\alpha_k^-}^0(2,d)$ and $\cN_{d-\alpha_k^+}(1,d-\alpha_k)\times\cN_{\alpha_k^-}(1,\alpha_k)$.
\endproof

\begin{remark}\label{Salpham-=S0alpham-}
 We have seen that in the case of $\cS_{\alpha_m^-}(2,d)$, the section $\theta_\gamma$ is always zero. Hence $\cS_{\alpha_m^-}(2,d)=\cS_{\alpha_m^-}^0(2,d)$ and therefore the similarity of Corollaries \ref{S-alpham} and \ref{S-alphak0}.
\end{remark}

\subsubsection{\texorpdfstring{$\cS_{\alpha_k^-}^1(2,d)$}{Salpha-1(2,d)}}We move on to the description of $\cS_{\alpha_k^-}^1(2,d)$. From Proposition \ref{q=0 implies ...}, $\cS_{\alpha_k^-}^1(2,d)$ is the space of pairs in $\cS_{\alpha_k^-}(2,d)$ such that extension (\ref{ext-}) is non-trivial.

Before going to the analysis of $\cS_{\alpha_k^-}^1(2,d)$ we first need the following proposition.

\begin{proposition}\label{S-rk2}
Let $\alpha_k<d/2$. If $(V,\gamma)\in\cS_{\alpha_k^-}(2,d)$, then generically $\rk(\gamma)=2$.
\end{proposition}
\proof
As always, $\gamma\neq 0$. Suppose that $\rk(\gamma)=1$. Then, $\det(\gamma)=0$, hence, from (\ref{det=bargamma2}), $\theta_\gamma=0$. From Proposition \ref{q=0 implies ...}, $V=M\oplus L$ and $$\gamma=\begin{pmatrix}
     \gamma' & 0 \\
     0 & \gamma''
\end{pmatrix}$$ with $\gamma''\neq 0$. Since also $\gamma'\neq 0$, it follows that $\det(\gamma)\neq 0$ and this is a contradiction with $\rk(\gamma)=1$.
\endproof

Given $((M,\gamma'),L)\in\cN_{d-\alpha_k^+}(1,d-\alpha_k)\times\Jac^{\alpha_k}(X)$ and recalling that $D=\divisor(\gamma')$,
consider the subvariety $$\mathcal C((M,\gamma'),L)$$ of $$H^0(D,M^{-1}L^{-1}U)\setminus\{0\}\times H^0(X,M^{-2}L^{-2}U^2)\setminus\{0\}$$ whose elements $(q,\eta)$ satisfy the equation $$q^2+\eta|_D=0.$$ 
$\C^*$ acts freely on $\mathcal C((M,\gamma'),L)$ as $$\lambda\cdot(q,\eta)=(\lambda q,\lambda^2\eta)$$ and we denote the quotient by
\begin{equation}\label{Q}
\mathcal Q((M,\gamma'),L)=\mathcal C((M,\gamma'),L)/\C^*.
\end{equation}

\begin{proposition}\label{S-alphak1}
Let $d_U$ even and $\alpha_k\neq\alpha_m,\alpha_M$ or $d_U$ odd and $\alpha_k\neq\alpha_M$. Suppose that $d_U-d>g-1$. Then there is a morphism $$\cS_{\alpha_k^-}^1(2,d)\longrightarrow\cN_{d-\alpha_k^+}(1,d-\alpha_k)\times\Jac^{\alpha_k}(X)$$ whose fibre over $((M,\gamma'),L)$ is isomorphic to $\mathcal Q((M,\gamma'),L)$ as defined in (\ref{Q}).
\end{proposition}
\proof
If $(V,\gamma)\in\cS_{\alpha_k^-}^1(2,d)$, we already know that we can write $V$ as the extension (\ref{ext-}), and that the pair $(M,\gamma')$ is $(d-\alpha_k^+)$-stable.
So we have the map
$$p:\cS_{\alpha_k^-}^1(2,d)\longrightarrow\cN_{d-\alpha_k^+}(1,d-\alpha_k)\times\Jac^{\alpha_k}(X)$$
given by $$p(V,\gamma)=((M,\gamma'),L).$$

Let $(V,\gamma)\in p^{-1}((M,\gamma'),L)$ and $D=\divisor(\gamma')$. Then $$\theta_\gamma\in H^0(D,M^{-1}L^{-1}U)\setminus\{0\}$$ because $(V,\gamma)\in\cS_{\alpha_k^-}^1(2,d)$ and, by the previous proposition, $$\det(\gamma)\in H^0(X,M^{-2}L^{-2}U^2)\setminus\{0\}.$$ Moreover,
$$\theta_\gamma^2+\det(\gamma)|_D=0$$ so we have the map 
\begin{equation}\label{fibre to set of square roots I1}
p^{-1}((M,\gamma'),L)\longrightarrow\mathcal Q((M,\gamma'),L)
\end{equation} given by
\begin{equation}\label{fibre to set of square roots I2}
 (V,\gamma)\mapsto [(\theta_\gamma,\det(\gamma))].
\end{equation}

Let us now see that we also have a map the other way around and which is inverse of the above one.
As we are assuming $d_U-d>g-1$, we have $\chi(M^{-2}L^{-2}U^2)>0$, hence $H^0(X,M^{-2}L^{-2}U^2)\neq 0$ and then $$\mathcal Q((M,\gamma'),L)\neq\emptyset.$$
Take $[(q,\eta)]\in\mathcal Q((M,\gamma'),L)$ and choose a representative $(q,\eta)$. We construct a pair $(V,\gamma)$ as follows.

Consider the following complexes 
$$C^\bullet_1:L^{-1}M\stackrel{\cong}{\longrightarrow}L^{-1}M^{-1}U(-D)$$
$$C^\bullet_2:L^{-1}M\stackrel{c}{\longrightarrow}L^{-1}M^{-1}U$$
and $$C^\bullet_3:0\longrightarrow L^{-1}M^{-1}U|_D$$
where 
\begin{equation}\label{c of complex}
c(\psi)=\gamma'\psi.
\end{equation}
We have the short exact sequence 
$$0\longrightarrow C^\bullet_1\longrightarrow C^\bullet_2\stackrel{r(D)}{\longrightarrow}C^\bullet_3\longrightarrow 0$$
given by the commutative diagram of sheaves of holomorphic sections
$$\xymatrix{0\ar[d]&0\ar[d]\\
L^{-1}M\ar[r]^(.35){\cong}\ar[d]_{=}&L^{-1}M^{-1}U(-D)\ar[d]^{c}&\\
L^{-1}M\ar[r]^(.4){c}\ar[d]&L^{-1}M^{-1}U\ar[d]^{r(D)}&\\
0\ar[r]^(.4){0}\ar[d]&L^{-1}M^{-1}U|_D\ar[d]&\\
0&0.}\vspace{0,5cm}$$
From this we obtain a long exact sequence in hypercohomology of the complexes
\begin{equation}\label{leshyp}
\begin{split}
0&\longrightarrow\mathbb{H}^0(X,C^\bullet_1)\longrightarrow
\mathbb{H}^0(X,C^\bullet_2)\stackrel{r(D)}{\longrightarrow}\mathbb{H}^0(X,C^\bullet_3)\longrightarrow\mathbb{H}^1(X,C^\bullet_1)\longrightarrow \mathbb{H}^1(X,C^\bullet_2)\\
&\stackrel{r(D)}{\longrightarrow}\mathbb{H}^1(X,C^\bullet_3)\longrightarrow\mathbb{H}^2(X,C^\bullet_1)\longrightarrow
\mathbb{H}^2(X,C^\bullet_2)\stackrel{r(D)}{\longrightarrow}\mathbb{H}^2(X,C^\bullet_3)\longrightarrow 0
\end{split}
\end{equation}
from which it follows that $r(D)$ yields a natural isomorphism
\begin{equation}\label{isom H0D hyper H1C}
\mathbb{H}^1(X,C^\bullet_2)\stackrel{r(D)}{\cong}\mathbb{H}^1(X,C^\bullet_3)\cong H^0(D,L^{-1}M^{-1}U).
\end{equation}

Consider the element $(0,q)\in\mathbb{H}^1(X,C^\bullet_3)$ and the corresponding class $r(D)^{-1}(0,q)\in \mathbb{H}^1(X,C^\bullet_2)$. With respect to some open covering $(U_a)_a$ of $X$, choose a representative $$(\lambda_{ab},\gamma_a'')$$ of the class $r(D)^{-1}(0,q)\in\mathbb{H}^1(X,C^\bullet_2)$. Recall then that $\gamma'\lambda_{ab}=\gamma_b''-\gamma_a''$.

Let $V$ be the vector bundle defined by taking on each open $U_a$ the direct sum
\begin{equation}\label{reconstruction of V 1}
M|_{U_a}\oplus L|_{U_a} 
\end{equation}
and gluing over $U_{ab}$ through
the map 
\begin{equation}\label{reconstruction of V 2}
f_{ab}=\begin{pmatrix}
 1_M & \lambda_{ab} \\
 0 & 1_L
\end{pmatrix}.
\end{equation}

Also over each open $U_a$, consider the section of $H^0(U_a,S^2(M\oplus L)^*\otimes U)$ given, with respect 
to the decomposition (\ref{reconstruction of V 1}),
by
\begin{equation}\label{reconstruction of gamma}
\gamma_a=\begin{pmatrix}
 \gamma' & \gamma_a'' \\
 \gamma_a''^t\otimes 1_U & (\gamma_a''^t\otimes 1_U)\gamma'^{-1}\gamma_a''+\eta\gamma'^{-1}
\end{pmatrix}.
\end{equation}
Observe that, since $r(D)(U_a)(\gamma_a'')=q|_{D\cap U_a}$ and $q^2+\eta|_D=0$, then $$(\gamma_a''^t\otimes 1_U)\gamma'^{-1}\gamma_a''+\eta\gamma'^{-1}$$ is defined over $D$.
One has $\gamma_b=f_{ab}^t\gamma_a f_{ab}$, so the collection of symmetric maps $(\gamma_a)_a$ yields a global symmetric map $\gamma:V\to V^*\otimes U$.

So, from $((M,\gamma'),L)\in\cN_{d-\alpha_k^+}(1,d-\alpha_k)\times\Jac^{\alpha_k}(X)$ and $(q,\eta)$ both non-zero and such that $q^2+\eta|_D=0$, we have built a $U$-quadratic pair $(V,\gamma)$ over $X$ such that $\det(\gamma)=\eta$, $\theta_\gamma=q$,  which lies in $\cS^1_{\alpha_k^-}(2,d)$ and which is mapped onto $((M,\gamma'),L)$ by the map $p$.

Suppose now that we had a different choice of the representative of the class $[(q,\eta)]\in\mathcal Q((M,\gamma'),L)$, say $(\beta q,\beta^2\eta)$ with $\beta\in\C^*$. From (\ref{isom H0D hyper H1C}), this pair defines a new class in $\mathbb{H}^1(X,C_2^\bullet)$ whose representative is $(\beta\lambda_{ab},\beta\gamma''_{ab})$.
The vector bundle $\tilde V$ constructed again from (\ref{reconstruction of V 1}) and gluing by
$$\begin{pmatrix}
 1_M & \beta\lambda_{ab} \\
 0 & 1_L
\end{pmatrix}=\tilde f_{ab}$$ is isomorphic to $V$ through the isomorphism $g:\tilde V\to V$ defined locally by $$g_a=\begin{pmatrix}
 1_M & 0\\
 0 & \beta
\end{pmatrix}$$ because $f_{ab}g_a=g_b\tilde f_{ab}$. Moreover, we consider the section
of $H^0(U_a,S^2(M\oplus L)^*\otimes U)$ given by
$$\begin{pmatrix}
 \gamma' & \beta\gamma_a'' \\
 \beta\gamma_a''^t\otimes 1_U & \beta^2((\gamma_a''^t\otimes 1_U)\gamma'^{-1}\gamma_a''+\eta\gamma'^{-1})
\end{pmatrix}=\tilde\gamma_a$$ and we have again $\tilde\gamma_b=\tilde f_{ab}^t\tilde\gamma_a\tilde f_{ab}$, so we have the pair $(\tilde V,\tilde\gamma)\in p^{-1}((M,\gamma'),L)$. Since $g_a^t\gamma_a g_a=\tilde\gamma_a$ the isomorphism $g$ is indeed an isomorphism between  $(\tilde V,\tilde\gamma)$ and $(V,\gamma)$.

In other words, we have a map
$$\mathcal Q((M,\gamma'),L)\longrightarrow p^{-1}((M,\gamma'),L)$$ defined by
$$[(q,\eta)]\mapsto\text{ isomorphism class of }(V,\gamma)\text{ defined by }(\ref{reconstruction of V 1}),\, (\ref{reconstruction of V 2})\text{ and }(\ref{reconstruction of gamma}).$$
Clearly, this map is inverse of that defined in (\ref{fibre to set of square roots I1}) and (\ref{fibre to set of square roots I2}) and this gives an isomorphism $\mathcal Q((M,\gamma'),L)\cong p^{-1}((M,\gamma'),L)$.
We have then seen that $p$ is surjective with fibre isomorphic to $\mathcal Q((M,\gamma'),L)$. It remains to check that $p$ is a morphism.

Let $p_X:\cN_{d-\alpha_k^+}(1,d-\alpha_k)\times X\to X$ be the projection. From Remark \ref{degenerated cases} and Theorem I of \cite{gomez-sols:2000}, there is a universal $p_X^*U$-quadratic pair $(\mathcal L_1,\pmb{\gamma}')$ over $\cN_{d-\alpha_k^+}(1,d-\alpha_k)\times X$. Consider also the Poincaré line bundle $\mathcal L_2$, of degree $\alpha_k$, over $\Jac^{\alpha_k}(X)\times X$.
Let $$\mathrm{pr}_{13}:\cN_{d-\alpha_k^+}(1,d-\alpha_k)\times\Jac^{\alpha_k}(X)\times X\longrightarrow\cN_{d-\alpha_k^+}(1,d-\alpha_k)\times X$$
$$\mathrm{pr}_{23}:\cN_{d-\alpha_k^+}(1,d-\alpha_k)\times\Jac^{\alpha_k}(X)\times X\longrightarrow\Jac^{\alpha_k}(X)\times X$$
and $$\mathrm{pr}_{12}:\cN_{d-\alpha_k^+}(1,d-\alpha_k)\times\Jac^{\alpha_k}(X)\times X\longrightarrow\cN_{d-\alpha_k^+}(1,d-\alpha_k)\times\Jac^{\alpha_k}(X)$$ be the projections.
Consider the following sheaves over $\cN_{d-\alpha_k^+}(1,d-\alpha_k)\times\Jac^{\alpha_k}(X)$: $$R^0\mathrm{pr}_{12*}(\mathrm{pr}_{13}^*\mathcal{L}_1^{-1}\mathrm{pr}_{23}^*\mathcal{L}_2^{-1}p_X^*U|_\Delta)$$ where $\Delta\subset\cSym^{d_U-2d+2\alpha_k}(X)\times X$ is the universal divisor, and $$R^0\mathrm{pr}_{12*}(\mathrm{pr}_{13}^*\mathcal{L}_1^{-2}\mathrm{pr}_{23}^*\mathcal{L}_2^{-2}p_X^*U^2).$$ Since we are assuming $d_U-d>g-1$, these spaces have constant dimension, hence are locally free. We consider the subsheaf $\widetilde{\mathcal C}((\mathcal L_1,\pmb{\gamma}'),\mathcal L_2)$ (of sets) of $$R^0\mathrm{pr}_{12*}(\mathrm{pr}_{13}^*\mathcal{L}_1^{-1}\mathrm{pr}_{23}^*\mathcal{L}_2^{-1}p_X^*U|_\Delta)\setminus\{0\}\times R^0\mathrm{pr}_{12*}(\mathrm{pr}_{13}^*\mathcal{L}_1^{-2}\mathrm{pr}_{23}^*\mathcal{L}_2^{-2}p_X^*U^2)\setminus\{0\}$$ consisting of pairs of non-zero sections
$(\mathbf q,\pmb{\eta})$
satisfying the equation $$\mathbf q^2+\pmb{\eta}|_\Delta=0.$$ If 
$\widetilde{\mathcal Q}((\mathcal L_1,\pmb{\gamma}'),\mathcal L_2)$ denotes the sheaf obtained from $\widetilde{\mathcal C}((\mathcal L_1,\pmb{\gamma}'),\mathcal L_2)$ by identifying sections of the form $(\mathbf q,\pmb{\eta})$ and $(\lambda\mathbf q,\lambda^2\pmb{\eta})$ for some $\lambda\in\C^*$, then this is a locally trivial fibration over $\cN_{d-\alpha_k^+}(1,d-\alpha_k)\times\Jac^{\alpha_k}(X)$ such that its fibre over $((M,\gamma'),L)$ is $\mathcal Q((M,\gamma'),L)$, as defined in (\ref{Q}). As in the proof of Proposition \ref{S+} (see also Proposition $3.4$ of \cite{thaddeus:1994}), we have the following commutative diagram:
$$\xymatrix{\widetilde{\mathcal Q}((\mathcal L_1,\pmb{\gamma}'),\mathcal L_2)\ar[r]^{\cong}\ar[rd]&\mathcal S_{\alpha_k^-}(2,d)\ar[d]^{p}\\
&\cN_{d-\alpha_k^+}(1,d-\alpha_k)\times\Jac^{\alpha_k}(X).}$$ So $p$ is a morphism and the result follows.
\endproof

One consequence of Corollary \ref{S-alphak0} and of the previous proposition is that $\cS_{\alpha_k^-}(2,d)=\cS^0_{\alpha_k^-}(2,d)\sqcup \cS^1_{\alpha_k^-}(2,d)$ is compact. It follows from Proposition \ref{S+- is a subvariety} (see also Remark \ref{S+- is a subvariety of N}) is a disconnected subvariety of $\cN_{\alpha_k^-}(2,d)$. Hence we can compute its dimension.

From Corollaries \ref{dimNalpha(1,d)} and \ref{S-alphak0}, we have
\begin{equation}\label{dimS0-}
\dim\cS_{\alpha_k^-}^0(2,d)=2d_U-2d.
\end{equation}

On the other hand, for $\cS^1_{\alpha_k^-}(2,d)$, we have:
\begin{corollary}
 If $d_U-d>g-1$, then $\dim\cS_{\alpha_k^-}^1(2,d)=3d_U-4d+2\alpha_k$.
\end{corollary}
\proof Since $d_U-d>g-1$, Proposition \ref{S-alphak1} holds.
From Corollary \ref{dimNalpha(1,d)} we have 
\begin{equation}\label{dim N(1,d-alpha)xJac}
\dim(\cN_{d-\alpha_k^+}(1,d-\alpha_k)\times\Jac^{\alpha_k}(X))=2\alpha_k-2d+d_U+g. 
\end{equation}

Given $((M,\gamma'),L)\in\cN_{d-\alpha_k^+}(1,d-\alpha_k)\times\Jac^{\alpha_k}(X)$,
 we now compute $$\dim\mathcal Q((M,\gamma'),L)=\dim\mathcal C((M,\gamma'),L)-1.$$ 
If $F:H^0(D,M^{-1}L^{-1}U)\setminus\{0\}\times H^0(X,M^{-2}L^{-2}U^2)\setminus\{0\}\to H^0(D,M^{-2}L^{-2}U^2)$
is given by $$F(q,\eta)=q^2+\eta|_D$$ then $\mathcal C((M,\gamma'),L)=F^{-1}(0)$. Linearising the map $F$ at a point $(q,\eta)$, we are lead to
$F_{(q,\eta)*}:H^0(D,M^{-1}L^{-1}U)\times H^0(X,M^{-2}L^{-2}U^2)\to H^0(D,M^{-2}L^{-2}U^2)$ with $$F_{(q,\eta)*}(\dot q,\dot\eta)=2q\dot q+\dot\eta|_D.$$
Choose $(q,\eta)$ such that $F_*$ is surjective (for instance, $(1,\eta)$, where $q=1$ means that $q(p)=1$ for each $p\in\cSupp(D)$).
Then $\dim\mathcal C((M,\gamma'),L)$ is
$$\dim\ker(F_{(1,\eta)*})=\dim H^0(X,M^{-2}L^{-2}U^2)=2(d_U-d)+1-g$$
because $d_U-d>g-1$ implies $\dim H^1(X,M^{-2}L^{-2}U^2)=0$. Hence
\begin{equation}\label{dim A((M,gamma'),L)}
\dim\mathcal Q((M,\gamma'),L)=2(d_U-d)-g.
\end{equation}
The result now follows from (\ref{dim N(1,d-alpha)xJac}) and (\ref{dim A((M,gamma'),L)}).
\endproof

From this, from (\ref{dimS0-}) and from $\alpha_k\geq\alpha_m=d-[d_U/2]$ follows that $$\dim\cS_{\alpha_k^-}^1(2,d)\geq\dim\cS_{\alpha_k^-}^0(2,d)$$ with equality if and only if $\alpha_k=\alpha_m$ and $d_U$ even. Hence from (\ref{dijunion}) and the previous corollary, we conclude the following.
\begin{corollary}\label{dimSalphak-(2,d)}
Let $\alpha\neq\alpha_M$ and $d$ such that $d_U-d>g-1$. Then each connected component of $\cS_{\alpha_k^-}(2,d)$ has dimension less or equal than $3d_U-4d+2\alpha_k$.
\end{corollary}

 \section{The space \texorpdfstring{$\cN_{\alpha_m^-}(2,d)$}{Nalpham-(2,d)}}\label{The space Nalpham-(2,d)}

Having examined the differences which occur on the moduli spaces when we cross a critical value of $\alpha$, we now address the problem of studying the number of connected components of one of them. This will be done in this section and the moduli space which will be analysed is $\cN_{\alpha_m^-}(2,d)$, the one for which the parameter $\alpha$ is less than the minimum critical value $\alpha_m=d-[d_U/2]$. In Section \ref{Number of connected components} we join the results of this and the previous sections to achieve the goal of computing the number of connected components of $\cN_\alpha(2,d)$, for any $\alpha\leq d/2$.

The method we shall employ to give the desired description of $\cN_{\alpha_m^-}(2,d)$ is the theory of spectral curves together with an analogue of the Hitchin map which is slightly outlined in the next sections.

\subsection{The spectral curve}\label{sec:spectral curve}

We shall give a rough description of the spectral curve of $X$ corresponding to a line bundle $L$ and a section of $L^2$. Then we shall see how to associate a spectral curve to 
a quadratic pair $(V,\gamma)$, with $\gamma$ generically non-degenerate.
The classical references for this theory, particularly its relations with Higgs bundles, are \cite{beauville-narasimhan-ramanan:1989,hitchin:1987 integrable systems}.

Let then $L$ be a holomorphic line bundle over $X$ with $\deg(L)>0$. We begin by reviewing the construction of the spectral curve $X_{s,L}$
associated to a section $s\in H^0(X,L^2)$. Consider the complex
surface $T$ given by the total space of the line bundle $L$, and let
$\pi:T\to X$ be the projection. The pullback $\pi^*L$ of $L$ to its
total space has a tautological section
$$\lambda\in H^0(T,\pi^*L)$$ defined by $\lambda(x)=x$.

\begin{definition}
  Let $s\in H^0(X,L^2)$. The \emph{spectral curve} $X_{s,L}$ associated to
  $s$ is the zero scheme in the surface $T$ of the section
  $$\lambda^2+\pi^*s\in H^0(T,\pi^*L^2).$$
\end{definition}

\begin{remark}
  \label{rem:singular-smooth-X_{s,L}}
  In the present case, the spectral curve $X_{s,L}$ is always reduced, but it may be singular
  and reducible. In fact, it is smooth if and only if $s$ only has
  simple zeros and it is irreducible if and only if $s$ is not the
  square of a section of $L$.
\end{remark}
\begin{remark}
  \label{rem:def-X_{s,L}}
  The above definition of spectral curve is a very particular case of a general definition. In fact, one can define a spectral curve associated to an element of the sum $\bigoplus_{k=1}^n H^0(X,L^k)$. See \cite{beauville-narasimhan-ramanan:1989,hitchin:1987 integrable systems}.
\end{remark}

\subsection{An analogue of the Hitchin map and its generic fibre}\label{The Hitchin map}
Consider a $U$-quadratic pair $(V,\gamma)\in\cN_{\alpha_m^-}(2,d)$. By Proposition \ref{injectivity}, $\det(\gamma)$ is a non-zero holomorphic section of $\Lambda^2V^{-2}U^2$ and one can consider its divisor $\divisor(\det(\gamma))\in \cSym^{2d_U-2d}(X)$.
Let $$\mathcal P_X$$ be the $2^{2g}$-cover of $\cSym^{2d_U-2d}(X)$ which fits in the commutative diagram
\begin{equation}\label{def: PX}
\xymatrix{&\mathcal P_X\ar[r]\ar[d]&\Jac^{d_U-d}(X)\ar[d]^{L\mapsto L^2}\\
&\cSym^{2d_U-2d}(X)\ar[r]_{D\mapsto\mathcal{O}(D)}&\Jac^{2d_U-2d}(X).}
\end{equation}
In other words, $$\mathcal P_X=\cSym^{2d_U-2d}(X)\times_{\Jac^{2d_U-2d}(X)}\Jac^{d_U-d}(X)$$ i.e., it is the fibred product of $\cSym^{2d_U-2d}(X)$ and $\Jac^{d_U-d}(X)$ over $\Jac^{2d_U-2d}(X)$, and its elements are pairs $(D,L)\in \cSym^{2d_U-2d}(X)\times\Jac^{d_U-d}(X)$ such that $$\mathcal{O}(D)\cong L^2.$$

In order to describe $\cN_{\alpha_m^-}(2,d)$, we shall use the following map, which is analogue to the so-called Hitchin map defined for the first time by Hitchin in \cite{hitchin:1987 integrable systems}.
Consider then the map
\begin{equation}\label{def of Hitchin map}
\begin{array}{rccl}
  h:&\cN_{\alpha_m^-}(2,d) & \longrightarrow & \mathcal P_X\\
  &(V,\gamma) & \longmapsto & (\divisor(\det(\gamma)),\Lambda^2V^{-1}U).
\end{array}
\end{equation}

\begin{definition}
  An \emph{$L$-twisted Higgs pair} of type
  $(n,d)$ over $X$ is a pair $(V,\varphi)$, where $V$ is a holomorphic
  vector bundle over $X$, with $\rk(V)=n$ and $\deg(V)=d$, and
  $\varphi$ is a global holomorphic section of $\End(V)\otimes L$,
  called the \emph{Higgs field}.

 Two $L$-twisted Higgs pairs $(V,\varphi)$ and
  $(V',\varphi')$ are \emph{isomorphic} if there is a holomorphic
  isomorphism $f:V\to V'$ such that $\varphi'f=(f\otimes 1_L)\varphi$.
\end{definition}

\begin{definition}\label{xi twisted associated to pair}
Let $(V,\gamma)$ be a $U$-quadratic pair of type $(2,d)$ over $X$ and let $\xi=\Lambda^2V^{-1}U$. The \emph{$\xi$-twisted Higgs pair $(V,\varphi)$ associated to $(V,\gamma)$} is the one induced from $(V,\gamma)$ and from the isomorphism
\begin{equation}\label{isomorphism Vdual VtensordetV}
g:V\otimes\xi\stackrel{\cong}{\longrightarrow} V^*\otimes U
\end{equation}
given by $$g(v\otimes\phi\otimes u)=\phi(v\wedge -)\otimes u,$$
where $v\otimes\phi\otimes u\in V\otimes\xi=V\otimes\Lambda^2V^{-1}U$. In other words, $\varphi=g^{-1}\gamma$.
\end{definition}

\begin{lemma}\label{comparison between isomorphism notions of quadratic pair and corresponding twisted Higgs bundles}
Suppose that $V$ and $V'$ are rank $2$ holomorphic vector bundles with the same determinant. Let $\xi=\Lambda^2V^{-1}U$. Let $(V,\gamma)$ and $(V',\gamma')$ be two $U$-quadratic pairs, and $(V,\varphi)$ and $(V',\varphi')$ be, respectively, the associated $\xi$-twisted Higgs, as in Definition \ref{xi twisted associated to pair}. 
\begin{enumerate}
 \item If $(V,\varphi)$ is isomorphic to $(V',\varphi')$ as $\xi$-twisted Higgs pairs, then $(V,\gamma)$ is isomorphic to $(V',\gamma')$ as $U$-quadratic pairs.
 \item If $(V,\gamma)$ is isomorphic to $(V',\gamma')$ as $U$-quadratic pairs, then there is some $\lambda\in\C^*$ such that $(V,\varphi)$ is isomorphic to $(V',\lambda\varphi')$ as $\xi$-twisted Higgs pairs.
\end{enumerate}
\end{lemma}
\proof
Let $f:V\to V'$ be an isomorphism between $(V,\varphi)$ and $(V',\varphi')$, that is, \begin{equation}\label{f isomo twisted Higgs bundles}
\varphi'f=(f\otimes 1_\xi)\varphi.
\end{equation}
Since $\Lambda^2V=\Lambda^2V'$, then $\det(f)=\lambda\in\C^*$.

Let $g:V\otimes\xi\to V^*\otimes U$ be the isomorphism (\ref{isomorphism Vdual VtensordetV}), and define $g':V'\otimes\xi\to V'^*\otimes U$ similarly.
Now, we have that
\begin{equation*}
\begin{split}
(f^t\otimes 1_U)g'(f\otimes 1_\xi)(v\otimes\phi\otimes u)&=(f^t\otimes 1_U)g'(f(v)\otimes\phi\otimes u)\\
&=(f^t\otimes 1_U)(\phi(f(v)\wedge -)\otimes u)\\
&=\phi(f(v)\wedge f(-))\otimes u\\
&=(\phi\det(f))(v\wedge -)\otimes u\\
&=\lambda\phi(v\wedge -)\otimes u
\end{split}
\end{equation*}
so we conclude, from the definition of $g$ in (\ref{isomorphism Vdual VtensordetV}), that
$$(f^t\otimes 1_U)g'(f\otimes 1_\xi)=\lambda g.$$
From this, from (\ref{f isomo twisted Higgs bundles}) and noticing that $g\varphi=\gamma$ and $g'\varphi'=\gamma'$, we conclude that $$(f^t\otimes 1_U)\gamma'f=\lambda\gamma.$$ Thus $\sqrt{\lambda^{-1}}f$ is an isomorphism between $(V,\gamma)$ and $(V',\gamma')$ and this settles the first item.

For the second item, if $f:V\to V'$ is an isomorphism between $(V,\gamma)$ and $(V',\gamma')$ then $(f^t\otimes 1_U)\gamma'f=\gamma$. It follows, as above, that $\varphi'f=\lambda^{-1}(f\otimes 1_\xi)\varphi$ where $\C^*\ni\lambda=\det(f)$. So $$(\lambda\varphi')f=(f\otimes 1_\xi)\varphi$$ and $f$ is an isomorphism between $(V,\varphi)$ and $(V',\lambda\varphi')$.
\endproof

\begin{definition}
A $\xi$-twisted Higgs pair $(V,\varphi)$ of type $(2,d)$ is \emph{semistable} if $\deg(L)\leq d/2$ for any line subbundle $L\subset V$ such that $\varphi(L)\subset L\xi$.
\end{definition}

\begin{proposition}\label{equivalence of stability quad pairs an Higgs pairs}
Let $(V,\gamma)$ be a $U$-quadratic pair of type $(2,d)$ and $\xi=\Lambda^2V^{-1}U$. Let $(V,\varphi)$ be the corresponding $\xi$-twisted Higgs pair, in the sense of Definition \ref{xi twisted associated to pair}. Then $(V,\gamma)$ is $\alpha_m^-$-semistable if and only if $(V,\varphi)$ is semistable.
\end{proposition}
\proof Assume that $(V,\varphi)$ is semistable and let $L\subset V$. As $(V,\gamma)$ is $\alpha_m^-$-semistable, then Proposition \ref{injectivity} says that $\gamma(L)\neq 0$.
Suppose that $\gamma(L)\subset L^\perp U$. It is easy to see that
\begin{equation}\label{eq:eq gamma varphi}
\gamma(L)\subset L^\perp U\Longleftrightarrow\varphi(L)\subset L\xi,
\end{equation}
and since $(V,\varphi)$ is semistable, it follows that $\deg(L)\leq d/2$.

Finally, suppose that $\gamma(L)\not\subset L^\perp U$, and $\deg(L)>d-\alpha_m^-$. Then $L$ is a destabilizing subbundle for $\alpha_m^-$. So by Proposition \ref{stabilization}, $(V,\gamma)$ is $\alpha$-unstable for every $\alpha<\alpha_m^-$ and, from above, the destabilizing subbundle must also be of type \textbf{(C)}. We see that for any $\alpha<\alpha_m$, there is $L'\subset V$ such that $\deg(L')>d-\alpha$. Letting $\alpha\to-\infty$ this contradicts the fact that the degrees of subbundles of $V$ are bounded above (see Corollary 10.9 of \cite{mukai:2003}).
We conclude that $(V,\gamma)$ is $\alpha_m^-$-semistable.

The proof of the other direction is straighforward, using (\ref{eq:eq gamma varphi}).
\endproof

If $\xi=\Lambda^2V^{-1}U$ and $(V,\varphi)$ is a $\xi$-twisted Higgs pair, consider the sections defined by the coefficients of the characteristic polynomial of $\varphi$: $$(-\tr(\varphi),\det(\varphi))\in H^0(X,\xi)\oplus H^0(X,\xi^2).$$ 
We have $\det(\varphi)=\det(\gamma)$ and, as $\gamma$ is symmetric, $\varphi$ has trace zero. Hence one can view $h(V,\gamma)=(\divisor(\det(\gamma)),\Lambda^2V^{-1}U)$ in (\ref{def of Hitchin map}) as given by $\xi$ and by the divisor of the section given by the characteristic polynomial of $\varphi$.
The spectral curve $X_{s,\xi}$ associated to $\xi$ and to the section $s=\det(\gamma)\in H^0(X,\xi^2)$ is
the curve inside the total space $T$ of $\xi$ defined by the equation
$$\lambda^2+\pi^*\det(\gamma)=0.$$

Now, let $(D,\xi)$ be any pair in $\mathcal P_X$, defined in (\ref{def: PX}). We want to describe the fibre of $h$ over $(D,\xi)$, i.e., the space of isomorphism classes of $\alpha_m^-$-semistable $U$-quadratic pairs $(V,\gamma)$ with $\divisor(\det(\gamma))=D$ and $\Lambda^2V$ isomorphic to $U\xi^{-1}$.

From $(D,\xi)$ we have a section $s\in H^0(X,\mathcal{O}(D))=H^0(X,\xi^2)$, defined up to a non-zero scalar, and one can construct the spectral curve associated to this section $s$. We denote this spectral curve by $$X_{D,\xi}$$ (in Remark \ref{independence of choice of section} below we give an explanation of this notation).

Given a line bundle $\xi$, let $$\cM_{\xi}^{U\xi^{-1}}$$ denote the moduli space of $\xi$-twisted Higgs pairs of rank two, with fixed determinant $U\xi^{-1}$ and with traceless Higgs field.
In \cite{gothen-oliveira:2010}, we carry out a study of the singular fibre of the Hitchin map $\mathcal H$ defined in $\cM_L^\Lambda$ for any $L$ with positive degree and any $\Lambda$:
$$\begin{array}{rccl}
  \mathcal H:&\cM_L^\Lambda & \longrightarrow & H^0(X,L^2)\\
  &(V,\varphi) & \longmapsto & \det(\varphi).
\end{array}$$

Moreover we have the following proposition, which is immediate from Lemma \ref{comparison between isomorphism notions of quadratic pair and corresponding twisted Higgs bundles} and Proposition \ref{equivalence of stability quad pairs an Higgs pairs}:

\begin{proposition}\label{fibres quad and Higgs isomorphic}
Let $(D,\xi)\in \mathcal P_X$. Then $h^{-1}(D,\xi)\in\cN_{\alpha_m^-}(2,d)$ is isomorphic to $\mathcal H^{-1}(s)\in\cM_{\xi}^{U\xi^{-1}}$, where $s\in H^0(X,\xi^2)$ is such that $\divisor(s)=D$.
\end{proposition}

\begin{remark}\label{independence of choice of section}
Recall that we made a choice of a section $s$ associated to the divisor $D$ and this choice induces a choice of the corresponding spectral curve, as explained in section \ref{sec:spectral curve}. However, the fibre of $h$ does not depend of this choice, due to Lemma \ref{comparison between isomorphism notions of quadratic pair and corresponding twisted Higgs bundles}. In fact, if we had a different choice $\lambda s$, for some $\lambda\in\C^*$, then we would be working on the spectral curve $X_{\lambda s,\xi}:x^2+\pi^*\lambda s=0$ and we would be working with $\xi$-Higgs pairs of the form $(V,\sqrt{\lambda}\varphi)$, where $(V,\varphi)$ is a $\xi$-Higgs pair coming from $X_{s,\xi}$. But, although these two $\xi$-Higgs pairs are not isomorphic, the corresponding $U$-quadratic pairs $(V,\gamma)$ and $(V,\sqrt{\lambda}\gamma)$ are isomorphic. This yields an isomorphism between the fibres of $h$ using $X_{s,\xi}$ and $X_{\lambda s,\xi}$.
This is the reason why we denote ``the'' spectral curve associated to $(D,\xi)$ by $X_{D,\xi}$.
\end{remark}

The study of $\mathcal H^{-1}(s)$ in \cite{gothen-oliveira:2010} is done by considering the cases where $X_{s,\xi}$ is smooth, singular and irreducible, and singular and reducible. The smooth case is the generic one, and it is well known that the fibre $\mathcal H^{-1}(s)$ in that case is a torsor for the Prym variety of the spectral curve (cf. \cite{beauville-narasimhan-ramanan:1989,hitchin:1987 integrable systems}). The case of singular and irreducible spectral curve is carried out by a careful study of the compactification of the Jacobian of the singular spectral curve associated to $s$, using the relation between this Jacobian with that of its desingularization. Finally, the study of $\mathcal H^{-1}(s)$ when $X_{s,\xi}$ is reducible is done by a direct analysis of the eigenbundles of $\varphi$.

Theorem 8.1 in \cite{gothen-oliveira:2010} and Proposition \ref{fibres quad and Higgs isomorphic} imply the following:
\begin{theorem}\label{fibre of h connected and dimension}
 Let $(D,\xi)\in\mathcal P_X$. Then the fibre of $h:\cN_{\alpha_m^-}(2,d)\to \mathcal P_X$ over $(D,\xi)$ is connected and has dimension $d_U-d+g-1$.
\end{theorem}

\section{Components of \texorpdfstring{$\cN_\alpha(2,d)$}{Nalpha(2,d)}}\label{Number of connected components}
From Theorem \ref{fibre of h connected and dimension} and from the fact that $\mathcal P_X$ is connected and $\dim\mathcal P_X=2d_U-2d$, one concludes the following:

\begin{theorem}\label{connectedness of extremal space}
For every $d<d_U$, the space $\cN_{\alpha_m^-}(2,d)$ is connected and has dimension $3(d_U-d)+g-1$.
\end{theorem}

Hence the dimension of $\cN_{\alpha_m^-}(2,d)$ is the expected dimension given in (\ref{dimension of moduli of pairs}).

Before stating our main result, we need one final lemma.
In the following all spaces are assumed to be second countable and
Hausdorff (and thus metrizable). Thus compactness is equivalent to
sequential compactness.
\begin{lemma}\label{connectedness properties}
  Let $N^{\pm}$ be compact spaces and let $S^{\pm} \subset
  N^{\pm}$ be proper closed subspaces. Assume that
  $\overline{(N^{\pm}\setminus S^{\pm})} = N^{\pm}$ and that there is
  a homeomorphism $N^{+}\setminus S^{+} \cong N^{-}\setminus
  S^{-}$. If $N^{-}$ and $S^{+}$ are connected, then so is $N^{+}$.
\end{lemma}

\begin{proof}
  Let $U^{\pm} = N^{\pm} \setminus S^\pm$. Then $U^{\pm}$ are
  non-compact.

  Suppose now that $N^+$ is not connected. Then, since the closure of
  a connected set is connected, $U^+$ is not
  connected. Let $N^+ = N^+_1 \cup N^+_2$ be a decomposition into
  disjoint non-empty closed subsets. Then $U^+_1 = U^+ \cap N^+_1$ and
  $U^+_2 = U^+ \cap N^+_2$ are disjoint non-empty open subsets of $U^+ \cong
  U^-$. By the connectedness of $N^{-}$, the
  intersections $\overline{U^{+}_i} \cap S^{-}$ are non-empty, where we are considering closures in $N^-$. As
  above, this implies that $U^{+}_i$ is non-compact for
  $i=1,2$. Considering now the closures in $N^+$, we
  have $\overline{U^{+}_i} = N^{+}_i$ and it follows that
  $N^{+}_i \cap S^{+}$ is non-empty for $i=1,2$. This shows that $S^+$
  is disconnected, a contradiction. 
\end{proof}

Now we reach our main result about the moduli of quadratic pairs.

\begin{theorem}\label{Nalpha(d,2) connected}
Let $d$ be such that $d_U-d>g-1$. For every $\alpha\leq d/2$, the moduli space $\cN_\alpha(2,d)$ is connected.
\end{theorem}
\proof
By Theorem \ref{connectedness of extremal space}, $\cN_\alpha(2,d)$ is connected, for every $\alpha<\alpha_m$.
We will see the flip loci described in Section \ref{Flip loci} have sufficient high codimension so that they do not affect the number of components of adjacent moduli spaces.


Again by Theorem \ref{connectedness of extremal space}, $\cN_{\alpha_m^-}(2,d)$ has dimension $3(d_U-d)+g-1$.

From Corollary \ref{Salphak+(2,d)} we have
$$\dim\cS_{\alpha_m^+}(2,d)=d_U-d+2g-2$$
hence, as $d_U-d>g-1$,
we have
\begin{equation}\label{dimSalpham+<dimNalpham-}
\dim\cS_{\alpha_m^+}(2,d)<\dim\cN_{\alpha_m^-}(2,d).
\end{equation}
On the other hand, from Corollary \ref{dimSalphak-(2,d)}
every point in $\cS_{\alpha_m^-}(2,d)$ is contained in a component whose dimension is less or equal than $$3d_U-4d+2\alpha_m=3d_U-2[d_U/2]-2d$$
hence,
\begin{equation}\label{dimSalpham-<dimNalpham-}
\dim\cS_{\alpha_m^-}(2,d)<\dim\cN_{\alpha_m^-}(2,d).
\end{equation}
Using (\ref{N--S-=N+-S+}), we conclude that $\dim\cN_{\alpha_m^+}(2,d)=\dim\cN_{\alpha_m^-}(2,d)=3(d_U-d)+g-1$.

Now, observe that (\ref{dimSalpham+<dimNalpham-}) and (\ref{dimSalpham-<dimNalpham-}) are valid for all critical value $\alpha_k<\alpha_M$ and not just $\alpha_m$. Hence we conclude that, for all $\alpha<\alpha_M$,
\begin{equation}\label{dimension of all moduli spaces}
\dim\cN_{\alpha}(2,d)=3(d_U-d)+g-1.
\end{equation}

So, for all $\alpha_k<\alpha_M=d/2$,
\begin{equation}\label{codim S+}
\codim\cS_{\alpha_k^+}(2,d)=2(d_U-d)-g+1>g-1\geq 1
\end{equation}
and, from Corollary \ref{dimSalphak-(2,d)}, every point in $\cS_{\alpha_k^-}(2,d)$ is contained in a component whose codimension is greater or equal than
\begin{equation}\label{codim S-}
d+g-1-2\alpha_k>g-1\geq 1.
\end{equation}

Recall that the flip loci measure the difference between two moduli spaces whose parameter lie on opposite sides of a critical value. From Theorem \ref{connectedness of extremal space}, (\ref{codim S+}) and (\ref{codim S-}), we see that the spaces $\cN_{\alpha_k^\pm}(2,d)$ and $\cS_{\alpha_k^\pm}$ satisfy the conditions of Lemma \ref{connectedness properties}. From this it follows that $\cN_\alpha(2,d)$ is connected for every generic $\alpha$.

If $\alpha_k\neq\alpha_M$ is a critical value, we have two obvious continuous maps 
$$\pi_{\pm}:\cN_{\alpha_k^{\pm}}(2,d)\longrightarrow\cN_{\alpha_k}(2,d).$$ 
From the definition of the flip loci $$\cN_{\alpha_k}(2,d)=\pi_-(\cN_{\alpha_k^-}(2,d))\cup\pi_+(\cN_{\alpha_k^+}(2,d)).$$ From above, $\pi_-(\cN_{\alpha_k^-}(2,d))\cap\pi_+(\cN_{\alpha_k^+}(2,d))$ is non-empty and the images of $\pi_{\pm}$ are connected.
The conclusion is that $\cN_{\alpha_k}(2,d)$ is also connected.
\endproof

\section{An application to surface group representations}
\label{sec:higgs-surface-groups}

\subsection{Higgs bundles}

Let $H\subset G$ be a maximal compact subgroup, and let $H^\C$ be the complexification of $H$. The Cartan decomposition,
$\lieg=\lieh\oplus\liem$, of $\lieg$, yields a decomposition
$\liegc=\liehc\oplus\liemc$ of the corresponding complexified Lie
algebra. Then $\liemc$ is a representation of $H^\C$ via the
\emph{isotropy representation}
\begin{equation}\label{complex isotropy representation}
\iota:H^\C\longrightarrow\Aut(\liemc)
\end{equation}
obtained by restricting the adjoint representation of $G^\C$ on
$\liegc$.  If $E_{H^\C}$ is a principal $H^{\C}$-bundle over $X$, we
denote by $E_{H^\C}(\liemc)=E\times_{H^{\C}}\liemc$ the vector bundle,
with fibre $\liemc$, associated to the isotropy representation.  Let
$K=T^*X^{1,0}$ be the canonical line bundle of $X$.

\begin{definition}\label{definition of Higgs bundle}
  A \emph{$G$-Higgs bundle} over a compact Riemann surface $X$ is a
  pair $(E_{H^\C},\varphi)$ where $E_{H^\C}$ is a principal
  holomorphic $H^\C$-bundle over $X$ and $\varphi$ is a global
  holomorphic section of $E_{H^\C}(\liemc)\otimes K$, called the
  \emph{Higgs field}.
\end{definition}

A $G$-Higgs bundle $(E_{H^\C},\varphi)$ is topologically classified by
the topological invariant of the corresponding $H^\C$-bundle
$E_{H^\C}$ and, as the maximal compact subgroup of $H^\C$ is $H$, the
topological classification of $G$-Higgs bundles is the same as the one
of $H$-principal bundles. Thus, whenever $G$ is connected, the topological class of a $G$-Higgs
bundle is given by an element in $H^2(X,\pi_1 H) \cong \pi_1 H$.

In \cite{garcia-gothen-mundet:2008}, a general notion of
(semi,poly)stability of $G$-Higgs bundles was developed, allowing for
proving a Hitchin--Kobayashi correspondence between polystable
$G$-Higgs bundles and solutions to certain gauge theoretic equations
known as Hitchin's equations. On the other hand, Schmitt
\cite{schmitt:2004,schmitt:2005,schmitt:2008} introduced stability
conditions for \emph{decorated bundles} and used these in his general
Geometric Invariant Theory construction of moduli spaces. In
particular, Schmitt's constructions give moduli of $G$-Higgs bundles
for the groups considered in this paper, and his stability conditions
coincide with the ones relevant for the Hitchin--Kobayashi
correspondence.  It should be noted that the stability conditions
depend on a parameter $\alpha \in \sqrt{-1}\lieh \cap \liez$, where
$\liez$ is the centre of $\liehc$.
We denote by
\begin{displaymath}
  \mathcal{M}_d(X,G)
\end{displaymath}
the moduli space of semistable (for the parameter value $\alpha=0$)
$G$-Higgs bundles with topological invariant $d \in \pi_1 H$. As
usual, the moduli space $\mathcal{M}_d(X,G)$ can also be viewed as
parametrizing isomorphism classes of polystable $G$-Higgs bundles.

\subsection{Higgs bundles for the adjoint form of the symplectic group}

Let $\Sp(2n,\R)$ be the real symplectic group of linear automorphisms
of $\R^{2n}$ which preserve the standard symplectic form. The centre
of $\Sp(2n,\R)$ is $Z(\Sp(2n,\R)) = \Z/2$ and we denote by
$\PSp(2n,\R) = \Sp(2n,\R) / (\Z/2)$ the projectivization of
$\Sp(2n,\R)$. A maximal compact subgroup of $\PSp(2n,\R)$ is isomorphic to
$\U(n)/(\Z/2)$, so the Cartan decomposition for
$\mathfrak{psp}(2n,\C)=\mathfrak{sp}(2n,\C)$ is given by
$\mathfrak{sp}(2n,\C)=\mathfrak{gl}(n,\C)\oplus\liemc$ where
\begin{equation}\label{mc for PSp}
\liemc=\left\{\left( \begin{array}{cc}
0 & B \\
C & 0
\end{array}\right)\st B,C\in\mathfrak{gl}(n,\C),\ B^T=B,\ C^T=C\right\}\cong S^2\C^n\oplus S^2(\C^n)^*. 
\end{equation}

Hence a \emph{$\PSp(2n,\R)$-Higgs bundle} over a compact Riemann
surface $X$ is a pair $(E,\varphi)$, where $E$ is a holomorphic
principal $\GL(n,\C)/(\Z/2)$-principal bundle and $\varphi$ is a
holomorphic global section of the vector bundle
$E\times_{\GL(n,\C)/(\Z/2)}(S^2\C^n\oplus S^2(\C^n)^*)\otimes K$.

We want to work with holomorphic vector bundles, so we shall use a
very similar procedure to the one taken in \cite{oliveira:2010} for
$G=\PGL(n,\R)$. Consider the group $\Sp(2n,\R)\times\U(1)$, the normal
subgroup $\{(I_n,1),(-I_n,-1)\}\cong\Z/2\vartriangleleft
\GL(n,\R)\times\U(1)$ and the corresponding quotient
group $$\Sp(2n,\R)\times_{\Z/2}\U(1)=(\Sp(2n,\R)\times\U(1))/(\Z/2).$$

\begin{notation}\label{notation enhanced}
We shall write 
\begin{align*}
  \ESp(2n,\R)&=\Sp(2n,\R)\times_{\Z/2}\U(1),\\
  \EU(n)&=\U(n)\times_{\Z/2}\U(1), \\
  \EGL(n,\C)&=\GL(n,\C)\times_{\Z/2}\C^*.
\end{align*}
The ``E'' stands for enhanced or extended.
\end{notation}

The complexification of the maximal compact subgroup
$\overline{H} = \EU(n) \subset \ESp(2n,\R)$ is
$\overline{H}^{\C}=\EGL(n,\C)$. Also,
$\mathfrak{\overline{g}}^{\C}=\mathfrak{\overline{h}}^\C\oplus\mathfrak{\overline{m}}^{\C}$
where $\mathfrak{\overline{g}}^\C=\mathfrak{sp}(2n,\C)\oplus\C$,
$\mathfrak{\overline{h}}^\C=\mathfrak{gl}(n,\C)\oplus\C$ and
$\overline\liem^\C=\liemc\oplus\{0\}\cong\liemc$, where $\liemc$ is
given by (\ref{mc for PSp}), so
\begin{equation}\label{mc for ESp}
\overline\liem^\C\cong\left\{\left( \begin{array}{cc}
0 & B \\
C & 0 
\end{array}\right)\st B,C\in\mathfrak{gl}(n,\C),\ B^T=B,\ C^T=C\right\}. 
\end{equation}

\begin{definition}
An \emph{$\ESp(2n,\R)$-Higgs bundle} over $X$ is a pair $(E,\varphi)$,
where $E$ is a holomorphic principal $\EGL(n,\C)$-bundle and $\varphi\in
H^0(X,E\times_{\EGL(n,\C)}\overline\liem^\C\otimes K)$, where $\overline\liem^\C$ is given by (\ref{mc for ESp}).
\end{definition}

Consider the actions of $\EGL(n,\C)$ on $\C^n$ and on $\C$ induced,
respectively, by the group homomorphisms
\begin{equation}\label{actCn}
\EGL(n,\C)\longrightarrow\GL(n,\C),\ \ \ \ [(w,\lambda)]\mapsto\lambda w 
\end{equation}
and
\begin{equation}\label{actC}
\EGL(n,\C)\longrightarrow\C^*,\ \ \ \ [(w,\lambda)]\mapsto\lambda^2.
\end{equation}
Note that together these two actions define an isomorphism
\begin{equation}
  \label{eq:EGL-iso}
  \begin{aligned}
    \EGL(n,\C) &\xrightarrow{\cong} \GL(n,\C)\times \C^* \\
    [(w,\lambda)]&\mapsto (\lambda w, \lambda^2 ).
  \end{aligned}
\end{equation}

We have the following description of an $\ESp(2n,\R)$-Higgs bundle in
terms of vector bundles:
\begin{proposition}\label{bundles from action of enhanced symplectic}
  Let $(\overline{E},\overline\varphi)$ be an $\ESp(2n,\R)$-Higgs
  bundle on $X$.  Through the actions (\ref{actCn}) and (\ref{actC})
  of $\EGL(n,\C)$ on $\C^n$ and on $\C$, associated to
  $(\overline{E},\overline\varphi)$ there is a quadruple
  $(V,L,\beta,\gamma)$, where $V$ is a rank $n$ holomorphic vector
  bundle, $L$ is a holomorphic line bundle and $(\beta,\gamma)\in
  H^0(X,(S^2V\otimes L^{-1}\oplus S^2V^*\otimes L)\otimes K)$.

  Moreover, two $\ESp(2n,\R)$-Higgs bundles
  $(\overline{E}_\nu,\overline\varphi_\nu)$, $\nu=1,2$, are isomorphic
  if and only if and only if the corresponding quadruples
  $(V_\nu,L_\nu,\beta_\nu,\gamma_\nu)$ are isomorphic, i.e., there are
  isomorphisms $V_1 \cong V_2$ and $L_1 \cong L_2$ intertwining
  $(\beta_1,\gamma_1)$ and $(\beta_2,\gamma_2)$.
\end{proposition}

\proof
From the actions (\ref{actCn}) and (\ref{actC}) we define, respectively,
the vector bundle $V=\overline{E}\times_{\EGL(n,\C)}\C^n$ and
the line bundle $L=\overline{E}\times_{\EGL(n,\C)}\C$.

Consider the representations $\sigma:\EGL(n,\C)\to\GL(S^2\C^n)$ and
$\sigma^*:\EGL(n,\C)\to\GL(S^2(\C^n)^*)$ given
by $$\sigma([w,\lambda])(B)=wBw^T$$
and $$\sigma^*([w,\lambda])(C)=(w^T)^{-1}Cw^{-1}.$$ If
$\overline\iota:\EGL(n,\C)\to\GL(\overline\liem^\C)$ is the isotropy
representation of $\EGL(n,\C)$ on $\overline\liem^\C$, then it is
clear that $\overline\iota([(w,\lambda)])(A)=\iota([w])(A)$, where
$\iota$ is the isotropy representation of $\GL(n,\C)/(\Z/2)$ in
$\liemc$. It is easy to see
that $$\overline\iota=\sigma\oplus\sigma^*$$ hence, taking into
account the actions (\ref{actCn}) and (\ref{actC}), from $\sigma$ we
obtain the vector bundle $S^2V\otimes L^{-1}$ and from $\sigma^*$ the
vector bundle $S^2V^*\otimes L$.  The Higgs field $\overline\varphi\in
H^0(X,\overline E\times_{\overline H^\C}\overline\liem^\C\otimes K)$
is therefore given, in terms of $V$ and $L$ by two sections:
$$\beta\in H^0(X,S^2V\otimes L^{-1}K)\hspace{1cm}
\text{and}\hspace{1cm}\gamma\in H^0(X,S^2V^*\otimes LK).$$ The
final statement about isomorphism of quadruples follows from the
isomorphism (\ref{eq:EGL-iso}).
\endproof

We shall slightly abuse notation and also call a quadruple
$(V,L,\beta,\gamma)$ as introduced in the preceding proposition an
$\ESp(2n,\R)$-Higgs bundles.

\begin{remark}
  \label{rem:esp-sp}
  An $\ESp(2n,\R)$-Higgs bundle $(V,L,\beta,\gamma)$ with $L=\mathcal
  O$ is the same thing as an $\Sp(2n,\R)$-Higgs bundle
  $(V,\beta,\gamma)$(cf.\ \cite{garcia-gothen-mundet:2008}).
\end{remark}

Projection on the first factor gives a homomorphism
\begin{displaymath}
  \ESp(2n,\R) \longrightarrow \PSp(2n,\R) 
\end{displaymath}
and so, to any $\ESp(2n,\R)$-Higgs bundle, we can naturally associate a
$\PSp(2n,\R)$-Higgs bundle. Note that this association is given by
extension of structure group in the principal bundles via
the map $\EGL(n,\C)\stackrel{p}{\longrightarrow}\GL(n,\C)/(\Z/2)$,
where $p([(w,\lambda)])=[w]$ and that the Higgs fields $\beta$ and
$\gamma$ are unchanged since the map $p$ intertwines the identity map
between the respective isotropy representations (\ref{mc for ESp}) and
(\ref{mc for PSp}). 

The following result is very similar to Proposition~5.4 of
\cite{oliveira:2010}, so we omit the proof.

\begin{proposition}\label{no obstruction to lift to enhanced symplectic}
Every $\PSp(2n,\R)$-Higgs bundle $(E,\varphi)$ on $X$ lifts to an
$\ESp(2n,\R)$-Higgs bundle $(\overline{E},\varphi)$.
\end{proposition}

\begin{proposition}\label{prop:iso-ESP-PSP}
  Two $\ESp(2n,\R)$-Higgs bundle $(V_\nu,L_\nu,\beta_\nu,\gamma_\nu)$,
  $\nu=1,2$ give rise to isomorphic $\PSp(2n,\R)$-Higgs bundles if and
  only if there is a line bundle $M$ on $X$ such that the
  $\ESp(2n,\R)$-Higgs bundles $(V_1,L_1,\beta_1,\gamma_1)$ and
  $(V_2\otimes M,L_2\otimes M^2,\beta_2,\gamma_2)$ are isomorphic.
\end{proposition}

\begin{proof}
  The correspondence with isomorphism of the underlying bundles is
  immediate from their definition. The complete statement including
  the Higgs fields follows because the Higgs fields are unchanged
  under the correspondence between $\ESp(2n,\R)$ and
  $\PSp(2n,\R)$-Higgs bundles.
\end{proof}

In view of Propositions~\ref{no obstruction to lift to enhanced
  symplectic} and \ref{prop:iso-ESP-PSP} we can now work
interchangably with either isomorphism classes of $\PSp(2n,\R)$-Higgs
bundles or with equivalence classes of $\ESp(2n,\R)$-Higgs bundles
under the equivalence relation introduced in the latter
Proposition. Thus we have the following immediate corollaries
(analogous, respectively, to Proposition~5.3 and Corollary~5.1 of
\cite{oliveira:2010}).

\begin{corollary}\label{deglambda01}
  Given a $\PSp(2n,\R)$-Higgs bundle $(E,\varphi)$, it is possible to
  choose a lift of $(E,\varphi)$ to an $\ESp(2n,\R)$-Higgs bundle
  $(V,L,\beta,\gamma)$ such that $L$ is trivial or $\deg(L)=1$.
\end{corollary}


\begin{corollary}\label{lift to Sp(2n,R)-Higgs}
  Let $(E,\varphi)$ be a $\PSp(2n,\R)$-Higgs bundle and
  $(V,L,\beta,\gamma)$ be an $\ESp(2n,\R)$-Higgs bundle which is a
  lift of $(E,\varphi)$. Then $(E,\varphi)$ lifts to an
  $\Sp(2n,\R)$-Higgs bundle if and only if $\deg(L)$ is even.
\end{corollary}

Next we give the topological classification of
\texorpdfstring{$\PSp(2n,\R)$}{PSp(2n,R)} and
\texorpdfstring{$\ESp(2n,\R)$}{ESp(2n,R)} bundles.
Restriction of the isomorphism (\ref{eq:EGL-iso}) gives an isomorphism
\begin{equation}
  \label{eq:EUn-iso}
    \begin{aligned}
    \epsilon\colon\EU(n) &\xrightarrow{\cong} \U(n)\times \U(1) \\
    [(w,\lambda)]&\mapsto (\lambda w, \lambda^2 ).
  \end{aligned}
\end{equation}
Hence (using the standard
identification $\pi_1\U(n) \cong \Z$)
\begin{equation}
  \label{eq:pi_1-EUn}
  \pi_1\EU(n) \cong \Z\x\Z.
\end{equation}

This means that $\ESp(2n,\R)$-Higgs
bundles are classified by a pair of integers and, thinking of an
$\ESp(2n,\R)$-Higgs bundle as a quadruple $(V,L,\beta,\gamma)$, we see
from Proposition~\ref{bundles from action of enhanced symplectic} that
this pair can be identified with
\begin{displaymath}
  (\deg(V),\deg(L)) \in \Z \x \Z.
\end{displaymath}
Using the identification (\ref{eq:EUn-iso}), the natural projection
$\EU(n,\R) \to \U(n)/(\Z/2)$ takes the form
\begin{displaymath}
  \begin{aligned}
    \U(n)\x\U(1) \xrightarrow{\epsilon}\EU(n) &\to \U(n)/(\Z/2), \\
    (g,\mu) &\mapsto \left[\frac{g}{\sqrt{\mu}}\right].
  \end{aligned}
\end{displaymath}
Thus we have a short exact sequence
\begin{displaymath}
  1 \to \U(1) \to \U(n)\x\U(1) \to \U(n)/(\Z/2) \to 1.
\end{displaymath}
Again using the standard
identification $\pi_1\U(n) \cong \Z$, the associated homotopy sequence gives
\begin{displaymath}
  1 \to \Z \xrightarrow{(n\cdot,2\cdot)} \Z\x\Z \to
  \pi_1\bigl(\U(n)/(\Z/2)\bigr) \to 1.
\end{displaymath}
It follows that, 
\begin{equation}
  \label{eq:pi_1Un-modZ2}
  \pi_1\bigl(\U(n)/(\Z/2)\bigr) \cong
  \begin{cases}
    \Z\times\Z/2 &\text{for $n$ even},\\
    \Z&\text{for $n$ odd}.
  \end{cases}
\end{equation}
In particular, the composition
\begin{math}
  \U(n) \into \EU(n) \xrightarrow{\epsilon} \U(n)\x\U(1) \to \U(n)/(\Z/2)
\end{math}
induces an isomorphism of $\pi_1\U(n)$ onto the $\Z$-factor in
$\pi_1\bigl(\U(n)/(\Z/2)\bigr)$.
For $n$ even and $(V,L)$ with $(\deg(V),\deg(L))=(d_1,d_2)$, the
projection $\pi_1\EU(n) \to \pi_1(\U(n)/(\Z/2))$ can then easily be
calculated to be
\begin{equation}\label{induced in pi1 to U(n)/Z_2 in terms of vector bundles}
(d_1,d_2)\mapsto (d_1-nd_2/2,d_2\mod 2),
\end{equation}
in terms of the identifications (\ref{eq:pi_1-EUn}) and
(\ref{eq:pi_1Un-modZ2}). 
Thus, from the point of view of $\PSp(2n,\R)$-Higgs bundles it is more
natural to work with the following invariants when $n$ is even:
\begin{equation}
  \label{eq:ab-invariants}
  (a,b)(V,L,\beta,\gamma) =  (\deg(V)-n\deg(L)/2,\deg(L)).
\end{equation}
In terms of these
\begin{displaymath}
  \deg(V) = a+nb/2 \quad\text{and}\quad \deg(L)=b.
\end{displaymath}
Notice that for a line bundle $F$ we have 
\begin{displaymath}
  (a,b)(V,L,\beta,\gamma) = (a,b)(V\otimes F,LF^2,\beta,\gamma)
\end{displaymath}
which is consistent with Proposition~\ref{prop:iso-ESP-PSP}. 

\begin{remark}
  \label{rem:obstruction-lifting}
  From either point of view, we see that the obstruction to lifting a
  $\PSp(2n,\R)$-Higgs bundle to an $\Sp(2n,\R)$-Higgs bundle is given
  by the invariant $d_2 = b$ (cf.\ Remark~\ref{rem:esp-sp} and
  Corollary~\ref{lift to Sp(2n,R)-Higgs}).
\end{remark}

\begin{remark}\label{SO(2,3)}
  We remark that $\PSp(4,\R)$ is isomorphic to the group
  $\mathrm{SO}_0(2,3)$ (i.e.\ the connected component of the identity
  of $\mathrm{SO}(2,3)$).  As explained in
  \cite{bradlow-garcia-prada-gothen:2005}, an
  $\mathrm{SO}_0(2,3)$-Higgs bundle is given by the the data
  $(W,Q_W,F,\beta,\gamma)$ where $F$ is a line bundle, $(W,Q_W)$ is a
  rank $3$ vector bundle equipped with a non-degenerate $F^2$-valued
  quadratic form, $\beta$ is a section of $\Hom(W,F)\otimes K$ and
  $\gamma$ a section of $\Hom(W,F^{-1})\otimes K$.  The objects are
  classified by two invariants, namely the degree of $F$ (which is
  actually the Toledo invariant: see (\ref{eq:ab-SO23}) below) and the
  second Stiefel-Whitney class $w_2(W,Q_W)\in\Z/2$.

  Generalizing the construction of a $\mathrm{SO}_0(2,3)$-Higgs bundle
  from an $\Sp(4,\R)$-Higgs bundle given in
  \cite[Section~3.3]{bradlow-garcia-prada-gothen:2005}, we can obtain
  a $\mathrm{SO}_0(2,3)$-Higgs bundle from an $\ESp(4,\R)$-Higgs
  bundle, as follows: If $(V,L,\beta,\gamma)$ is an $\ESp(4,\R)$-Higgs
  bundle, define
  $$W=S^2V\otimes\Lambda^2V^{-1}$$ and $$F=\Lambda^2V\otimes L^{-1}.$$
  Then, if $Q(x\otimes y,x'\otimes y')=(x\wedge x')\otimes (y\wedge
  y')$, then $Q$ is an $F^2L^2$-valued quadratic form on $S^2V$, hence
  $W$ has the induced non-degenerate quadratic form $Q_W$. Moreover,
  $\gamma$ is a section of $\Hom(W,F^{-1})\otimes K$ and, since
  $W\cong W^*$, we can view $\beta$ as a section of $\Hom(W,F)\otimes
  K$ and. Hence we have obtained the $\mathrm{SO}_0(2,3)$-Higgs bundle
  $(W,Q_W,F,\beta,\gamma)$. It is easily checked that the invariants
  $(a,b)$ of the $\ESp(4,\R)$-Higgs bundle $(V,L,\beta,\gamma)$ are
  given by
  \begin{equation}
    \label{eq:ab-SO23}
    \begin{aligned}
      a &= \deg(L), \\
      b &= \deg(L) + w_2(W,Q_W) \mod 2.
    \end{aligned}
  \end{equation}
\end{remark}

\subsection{Stability, moduli spaces and the non-abelian Hodge Theorem}\label{stability}
 
In \cite{garcia-gothen-mundet:2008}, a general notion of
(semi,poly)stability for $G$-Higgs bundles was introduced and a
Hitchin--Kobayashi correspondence was established showing that
polystability of a $G$-Higgs bundles is equivalent to the existence of
a solution to certain gauge theoretic equations, known as the Hitchin
equations.  The general definition of stability is fairly involved but
in many examples it can be significantly simplified.  In the case of
$G=\ESp(2n,\R)$ a simplification can be carried out in a manner
entirely analogous to the case of $G=\Sp(2n,\R)$ studied in
\cite[Section~4]{garcia-gothen-mundet:2008} and the stability
condition then takes the following form.

\begin{proposition}\label{ssEsp2n}
An $\ESp(2n,\R)$-Higgs bundle $(V,L,\beta,\gamma)$ is 
\begin{itemize}
 \item semistable if and only if for any filtration of holomorphic subbundles $0\subset V_1\subset V_2\subset V$ such that
$$(\beta,\gamma)\in H^0(X,(S^2V_2+V_1\otimes_S V)\otimes L^{-1}K\oplus(S^2V_1^\perp+V_2^\perp\otimes_S V^*)\otimes LK),$$
we have $$\deg(V_1)+\deg(V_2)\leq\deg(V).$$
 \item stable if and only if for any filtration of holomorphic subbundles $0\subset V_1\subset V_2\subset V$ such that
$$(\beta,\gamma)\in H^0(X,(S^2V_2+V_1\otimes_S V)\otimes L^{-1}K\oplus(S^2V_1^\perp+V_2^\perp\otimes_S V^*)\otimes LK),$$
the following holds: if at least one of the subbundles $V_1$ or $V_2$ is proper, then $$\deg(V_1)+\deg(V_2)<\deg(V)$$ and in any other case,
$$\deg(V_1)+\deg(V_2)\leq\deg(V).$$
\end{itemize} 
\end{proposition}

\begin{remark}\mbox{}
\begin{enumerate}
\item The general notion of semistability of $G$-Higgs bundles depends
  on a parameter $\alpha\in \sqrt{-1}\lieh\cap\liez$ where $\liez$ is
  the centre of $\liehc$. For $G=\ESp(2n,\R)$, we have
  $\sqrt{-1}\lieh\cap\liez =
  \sqrt{-1}(\mathfrak{u}(1)\oplus\mathfrak{u}(1))=\R\times\R$, so
  there is a $\alpha_1,\alpha_2)$-semistability condition. This is
  very similar to the case of $G=\Sp(2n,\R)$. However, if
  $(V,L,\beta,\gamma)$ is an $\ESp(2n,\R)$-Higgs bundle, it can be
  seen that if $\alpha_2\neq\deg(L)$, then it is
  $(\alpha_1,\alpha_2)$-unstable. Therefore we are considering
  $\alpha_2=\deg(L)$ fixed. Also, in the preceding theorem, we have
  restricted ourselves to the case of $\alpha_1=0$ semistability,
  because it is for this value of the parameter that the fundamental
  correspondence between moduli spaces $G$-Higgs bundles and
  $G$-character varieties (Theorem~\ref{thm:non-ab-Hodge} below)
  holds.
\item The notion polystablity of $\ESp(2n,\R)$-Higgs bundles is
  analogous to the one in Proposition~4.16 of
  \cite{garcia-gothen-mundet:2008}, with the obvious modifications.
\end{enumerate}
\end{remark}

With regard to the relation between the stability conditions for
$\PSp(2n,\R)$-Higgs bundles and $\ESp(2n,\R)$-Higgs bundles we have
the following result.

\begin{proposition}
  \label{prop:ESp-PSp-stability-equivalence}
  An $\ESp(2n,\R)$-Higgs bundle $(V,L,\beta,\gamma)$ is polystable if
  and only if the associated $\PSp(2n,\R)$-Higgs bundle is polystable.
\end{proposition}

\begin{proof}
  This can be checked by specializing and comparing directly the
  general polystability conditions given in
  \cite{garcia-gothen-mundet:2008}. An alternative proof can be given
  by invoking the Hitchin--Kobayashi correspondence proved in that paper,
  since the existence of solutions to the Hitchin equations on an
  $\ESp(2n,\R)$-Higgs bundle is clearly equivalent to the existence of
  solutions on the corresponding $\PSp(2n,\R)$-Higgs bundle.
\end{proof}

Next we recall the non-abelian Hodge theory correspondence. Let $G$
be a connected semisimple real Lie group with maximal compact subgroup
$H \subseteq G$. By a \emph{representation} of $\pi_1
X$ in $G$ we mean a homomorphism $\rho\colon \pi_1 X \to G$. A
representation $\rho$ is \emph{reductive} if its composition with the
adjoint representation of $G$ on $\lie{g}$ is a completely reducible
representation. The \emph{character variety} for representations of
$\pi_1 X$ in $G$ is
\begin{displaymath}
  \mathcal{R}(\pi_1X,G) = \Hom^{\mathrm{red}}(\pi_1 X,G) /G,
\end{displaymath}
where $G$ acts by overall conjugation on homomorphisms. Any
representation $\rho$ has a topological invariant $c(\rho) \in \pi_1
H$ defined as the topological class of the associated flat bundle. 
  Let
  \begin{displaymath}
  \mathcal{R}_d(\pi_1X,G)\subset \mathcal{R}(\pi_1X,G)
\end{displaymath}
be the subspace of equivalence classes of
represensentations whose topological invariant is $c(\rho) =d$.

The non-abelian Hodge Theorem
(\cite{corlette:1988,donaldson:1987,hitchin:1987,simpson:1992,garcia-gothen-mundet:2008})
now states the following.

\begin{theorem}\label{thm:non-ab-Hodge}
  There is a homeomorphism
  \begin{displaymath}
    \mathcal{M}_d(X,G) \cong \mathcal{R}_d(\pi_1X,G).
  \end{displaymath}
\end{theorem}

\subsection{Bounds on invariants}

From now on we restrict to the case of $n$ being even.
Let $$\cM_{d_1,d_2}=\cM_{d_1,d_2}(X,\ESp(2n,\R))$$ the moduli space of
polystable $\ESp(2n,\R)$-Higgs bundles $(V,L,\beta,\gamma)$ with
$\deg(V)=d_1$ and $\deg(L)=d_2$.  Let also
$$\widehat\cM_{a,b}=\cM_{a,b}(X,\PSp(2n,\R)),$$ 
the moduli space of
polystable $\PSp(2n,\R)$-Higgs bundles with with topological invariants
$(a,b)\in\Z\times\Z/2$.

For a semistable $\ESp(2n,\R)$-Higgs bundle $(V,L,\beta,\gamma)$ with
$\deg(V)=d_1$ and $\deg(L)=d_2$, we have a Milnor-Wood inequality
(a Higgs bundle proof of this inequality can be easily given, cf.\
\cite{bradlow-garcia-prada-gothen:2003} for the case $G=\U(p,q)$ which
implies the result in the current setting):
\begin{equation}\label{bound on deg(V) 1}
n(1-g)+d_2\leq d_1\leq n(g-1)+d_2.
\end{equation}
This is  equivalent to
\begin{displaymath}
\abs{a} \leq n(g-1)
\end{displaymath} where $a=d_1-n d_2/2\in\Z$ is the invariant introduced in
(\ref{eq:ab-invariants}). In the context of surface group
representations the invariant $a$ is  the \emph{Toledo invariant}.

Furthermore, for $(d_1,d_2)\in\Z\times\Z$ such that (\ref{bound on
  deg(V) 1}) holds, we have the isomorphism 
\begin{displaymath}
\cM_{d_1,d_2}\cong\cM_{nd_2-d_1,d_2} 
\end{displaymath}
given by $(V,L,\beta,\gamma)\mapsto(V^*\otimes L,L,\gamma^t\otimes
1_K,\beta^t\otimes 1_K)$ and the induced isomorphism
\begin{displaymath}
  \widehat\cM_{a,b}\cong\widehat\cM_{-a,b}.
\end{displaymath}
We can, therefore, assume that
\begin{equation}\label{bound on deg(V)}
0 \leq a \leq n(g-1).
\end{equation}

Let now $L_0$ be a fixed line bundle of degree $1$ over $X$. Denote by $$\cM_{d,L_0}\subset\cM_{d,1}$$ be the subspace of $\ESp(4,\R)$-Higgs bundles $(V,L,\beta,\gamma)$ with $L=L_0$.
Similarly, let $$\cM_{d,\mathcal O}\subset\cM_{d,0}$$ be the subspace of $\ESp(4,\R)$-Higgs bundles $(V,L,\beta,\gamma)$ with $L$ isomorphic to the trivial line bundle $\mathcal O$.

From Corollary \ref{deglambda01} and (\ref{induced in pi1 to U(n)/Z_2
  in terms of vector bundles}) the following is clear:
\begin{proposition}\label{projmoduli}
  Let $d$ and $d'$ be two integers satisfying $0\leq d\leq n(g-1)$ and
  $1\leq d'\leq n(g-1)+1$. Let $[(V,L,\beta,\gamma)]$ denote the class
  of the corresponding $\PSp(4,\R)$-Higgs bundle under the equivalence
  relation given in Proposition \ref{prop:iso-ESP-PSP}.  Then the
  projection $(V,L,\beta,\gamma)\mapsto[(V,L,\beta,\gamma)]$ yields a
  continuous surjective map $$\cM_{d,\mathcal
    O}\sqcup\cM_{d',L_0}\longrightarrow\widehat\cM_{d,0}\sqcup\widehat\cM_{d'-1,1}$$
  preserving the decompositions.
\end{proposition}

\subsection{Relation with quadratic pairs and connectedness theorems}
For the remainder of the paper we specialize to the case $n=2$, i.e.,
$G=\ESp(4,\R)$ or $G=\PSp(4,\R)$.  Our goal is to count the number of
connected components of $\widehat\cM_{a,b}$ for $(a,b)\in\Z\times\Z/2$
such that $0 < \abs{a} < 2g-2$. The situation for $\abs{a}=0$ and
$\abs{a}=2g-2$ is somewhat special and, at any rate, in these cases the
count follows from the results of \cite{gothen:2001} and
\cite{bradlow-garcia-prada-gothen:2005}.  Note also that the count for
$b=0$ (corresponding to $\PSp(4,\R)$-Higgs bundles which lift to
$\Sp(4,\R)$-Higgs bundles) follows from the results of
Garc\'\i{}a-Prada and Mundet \cite{garcia-prada-mundet:2004}.

We will analyze the spaces $\cM_{d,0}$ and $\cM_{d,L_0}$ and from that draw our conclusions about $\widehat\cM_{a,b}$, using Proposition \ref{projmoduli}.
Let us deal first with $\cM_{d,L_0}$, with 
$$1<d< 2g-1.$$
We introduce the following \emph{Hitchin proper functional}. It is defined as
\begin{align*}
f:\cM_{d,L_0}&\longrightarrow\hspace{0,2cm}\R\\
(V,L_0,\beta,\gamma)\hspace{0,25cm}&\longmapsto \|\beta\|^2_{L^2}+\|\gamma\|^2_{L^2}
\end{align*}
\begin{remark} The definition of the Hitchin
  functional uses a \emph{harmonic metric} on $V$ coming from the
  Hitchin-Kobayashi correspondence --- see \cite{hitchin:1987}.
\end{remark}

The following consequence of properness of the non-negative function
$f$ is well known (cf.\ \cite{hitchin:1987} or Proposition~4.3 of
\cite{bradlow-garcia-prada-gothen:2003}).

\begin{proposition}\label{rem:local-minima}
  The space $\cM_{d,L_0}$ is connected if the subspace of local minima
  of the Hitchin proper function is connected.
\end{proposition}

Using very similar methods to the ones of
\cite{gothen:2001,garcia-prada-mundet:2004,garcia-gothen-mundet:2008
  II} for $G=\Sp(2n,\R)$, one can prove the following result:

\begin{proposition}
  Let $(V,L_0,\beta,\gamma)$ represent a point in $\cM_{d,L_0}$, with
  $1<d<2g-1$. Then $(V,L_0,\beta,\gamma)$ is a minimum of $f$ if and
  only if $\beta=0$.
\end{proposition}

We have the following immediate corollary.

\begin{proposition}\label{prop:M-local-minima}
    For any integer $1<d<2g-1$, the subvariety of local minima of $f$ is
  the moduli space $$\cN_{d,L_0}$$ of semistable $\ESp(4,\R)$-Higgs
  bundles $(V,L_0,0,\gamma)$ such that $V$ is a rank $2$ holomorphic
  vector bundle of degree $d$ and $\gamma\in H^0(X,S^2V^*\otimes
  L_0K)$.
\end{proposition}

Now, we connect this with the study of quadratic pairs made in the
first part of the paper. Let $$\cN_0(2,d)$$ be the moduli space of
$0$-semistable $L_0K$-quadratic pairs of type $(2,d)$.
\begin{proposition}\label{minimum = quadratic pairs}
The spaces $\cN_{d,L_0}$ and $\cN_0(2,d)$ are isomorphic.
\end{proposition}
\proof In view of Proposition~\ref{prop:M-local-minima} the result
follows by comparing the notions of $0$-(semi,poly)stability given in
Proposition \ref{sspairs2} for quadratic pairs, and from the notion of
(semi,poly)stability for $\ESp(4,\R)$-Higgs bundles in Theorem
\ref{ssEsp2n}. By considering all possible filtrations $0\subset
V_1\subset V_2\subset V$ of the rank $2$ bundle $V$ in Theorem
\ref{ssEsp2n}, one easily checks that these notions coincide.
\endproof

\begin{proposition}\label{connectedness of part deg 1}
 For each integer $d$ such that $3-2g<d<2g-1$ and $d\neq 1$, the space $\cM_{d,L_0}$ is connected.
\end{proposition}
\proof Recall that we can assume $1<d< 2g-1$.  By Theorem
\ref{Nalpha(d,2) connected} one has that $\cN_0(2,d)$ is connected for
every $1<d<g$, hence, by Proposition~\ref{minimum = quadratic pairs},
the same is valid for $\cN_{d,L_0}$ for such $d$.

If $g\leq d< 2g-1$, then $\cN_0(2,d)$ corresponds to the case
$\cN_{\alpha_m^-}(2,d)$, because in this case the formula for
$\alpha_m$ given in Notation \ref{formula alphak}, yields
$\alpha_m=d-g+1>0$. Hence, from Theorem \ref{connectedness of extremal
  space}, $\cN_0(2,d)$ is connected, so Proposition \ref{minimum =
  quadratic pairs}, says that $\cN_{d,L_0}$ is connected as well.
	
Now the result follows from Proposition~\ref{rem:local-minima}.\endproof

The connectedness of $\cM_{d,\mathcal{O}}$ with $0 < d < 2g-2$ was
proved by Garc\'\i{}a-Prada and Mundet (alternatively the argument
used above to prove Proposition~\ref{connectedness of part deg 1}
could be applied to give a proof):

\begin{proposition}[\protect{\cite[Theorem~5]{garcia-prada-mundet:2004}}]
   \label{connectedness of part deg 0}
  For each integer $d$ such that $0<|d|< 2g-2$, $\cM_{d,\mathcal{O}}$
  is connected.
\end{proposition}

We are now ready to state the theorem on the connectedness of the moduli space of $\PSp(4,\R)$-Higgs bundles, with fixed topological classes.

\begin{theorem}\label{connected components of cRPSp(4,R)}
  For each $(a,b)\in\Z\times\Z/2$ such that $0 < \abs{a}<2g-2$, the space
  $\widehat\cM_{a,b}$ is connected.
\end{theorem}
\proof 
Follows from Propositions \ref{connectedness of part deg 1},
\ref{connectedness of part deg 0} and \ref{projmoduli}.
\endproof

Using the non-abelian Hodge theory correspondence of
Theorem~\ref{thm:non-ab-Hodge}, we can rephrase our
Theorem~\ref{connected components of cRPSp(4,R)} as follows:
\begin{theorem}\label{connected-components-RPSp}
  For each $(a,b)\in\Z\times\Z/2$ such that $0<|a|<2g-2$, the space
  $\mathcal{R}_{a,b}(\pi_1X,\PSp(4,\R))$ is connected.
\end{theorem}

Recalling the correspondence of Remark~\ref{SO(2,3)}, we can
alternatively consider the character variety
\begin{displaymath}
  \mathcal{R}_{a,w}(\pi_1X,\mathrm{SO}_0(2,3))
\end{displaymath}
of representations of $\pi_1 X$ in $\mathrm{SO}_0(2,3)$ with
invariants $(a,w) \in \Z\times\Z/2$. We then have an identification
\begin{displaymath}
  \mathcal{R}_{a,w}(\pi_1X,\mathrm{SO}_0(2,3))
  = \mathcal{R}_{a,b}(\pi_1X,\PSp(4,\R))
\end{displaymath}
where the invariants are related by $(a,b) = (a,a+w \mod 2)$ (see
(\ref{eq:ab-SO23})). (There is of course an analogous identification of
the corresponding Higgs bundle spaces.) We thus have the following
equivalent formulation of Theorem~\ref{connected-components-RPSp}:

\begin{theorem}\label{connected-components-RSO23}
  For each $(a,w)\in\Z\times\Z/2$ such that $0<\abs{a}<2g-2$, the space
  $\mathcal{R}_{a,w}(\mathrm{SO}_0(2,3))$ is connected.
\end{theorem}


\begin{thebibliography}{99}


\bibitem{beauville-narasimhan-ramanan:1989}
{A. Beauville, M. S. Narasimhan, S. Ramanan}, \emph{Spectral curves and the generalised theta divisor}, J. Reine Angew. Math. \textbf{398} (1989), 169--179.

\bibitem{biswas-ramanan:1994}
{I. Biswas, S. Ramanan}, \emph{An infinitesimal study of the
moduli of Hitchin pairs}, J. London Math. Soc. (2) \textbf{49}
(1994), 219--231.

\bibitem{bradlow:1991}
S.~B. Bradlow, \emph{Special metrics and stability for holomorphic bundles with
  global sections}, J. Differential Geom. \textbf{33} (1991), 169--213.

\bibitem{bradlow-et-al:1995}
S.~Bradlow, G.~D. Daskalopoulos, O.~Garc{\'{\i}}a-Prada, and
  R.~Wentworth, \emph{Stable augmented bundles over {R}iemann surfaces},
  Vector bundles in algebraic geometry ({D}urham, 1993), London Math. Soc.
  Lecture Note Ser., vol. 208, Cambridge Univ. Press, Cambridge, 1995,
  pp.~15--67.

\bibitem{bradlow-garcia-prada-gothen:2003}
{S. B. Bradlow, O. Garc\'{\i}a-Prada, P. B. Gothen}, \emph{Surface
group representations and $\U(p,q)$-Higgs bundles}, J. Diff. Geom.
\textbf{64} (2003), 111--170.


\bibitem{bradlow-garcia-prada-gothen:2004 triples}
{S. B. Bradlow, O. Garc\'{\i}a-Prada, P. B. Gothen}, \emph{Moduli spaces of holomorphic triples over compact Riemann surfaces}, Math. Ann. \textbf{328} (2004),
299--351.

\bibitem{bradlow-garcia-prada-gothen:2005}
{S. B. Bradlow, O. Garc\'{\i}a-Prada, P. B. Gothen}, \emph{Maximal
surface group representations in isometry groups of classical
Hermitian symmetric spaces}, Geometriae Dedicata \textbf{122} (2006), 185--213.


\bibitem{bradlow-garcia-prada-mundet:2003}
{S. B. Bradlow, O. Garc\'{\i}a-Prada, I. Mundet i Riera},
\emph{Relative Hitchin-Kobayashi correspondences for principal
pairs}, Quart. J. Math. \textbf{54} (2003), 171--208.

\bibitem{corlette:1988}
{K. Corlette}, \emph{Flat $G$-bundles with canonical metrics},
J. Diff. Geom. \textbf{28} (1988), 361--382.

\bibitem{domic-toledo:1987}
{A. Domic, D. Toledo}, \emph{The Gromov norm of the Kaehler class of
  symmetric domains}, Math. Ann. \textbf{276} (1987), 425--432.

\bibitem{donaldson:1987}
{S. K. Donaldson}, \emph{Twisted harmonic maps and self-duality
equations}, Proc. London Math. Soc. (3) \textbf{55} (1987),
127--131.


\bibitem{garcia-gothen-mundet:2008}
{O. García-Prada, P. B. Gothen, I. Mundet i Riera}, The Hitchin-Kobayashi correspondence, Higgs pairs and surface group representations, Preprint arXiv:0909.4487v2.

\bibitem{garcia-gothen-mundet:2008 II}
{O. García-Prada, P. B. Gothen, I. Mundet i Riera}, Higgs bundles and surface group representations in the real symplectic group, Preprint arXiv:0809.0576v3.

 \bibitem{garcia-prada-gothen-munoz:2007}
 {O. García-Prada, P. B. Gothen, V. Muñoz}, \emph{Betti numbers of the moduli space of rank $3$ parabolic Higgs bundles}, Memoirs Amer. Math. Soc. \textbf{187} (2007).

\bibitem{garcia-prada-mundet:2004}
{O. Garc\'{\i}a-Prada, I. Mundet i Riera}, \emph{Representations of the fundamental group of a closed oriented surface in $\cSp(4,\R)$}, Topology \textbf{43} (2004), 831--855.


\bibitem{gomez-sols:2000}
{T. Gómez, I. Sols}, \emph{Stability of conic bundles}, Internat. J. Math. \textbf{11} (2000), 1027--1055.

\bibitem{gothen:2001}
{P. B. Gothen}, \emph{Components of spaces of representations and stable triples}, 
Topology \textbf{40} (2001), 823-850.

\bibitem{gothen-oliveira:2010}
{P. B. Gothen, A. G. Oliveira}, \emph{The singular fibre of the Hitchin map}, Preprint arXiv:1012.5541v2.


 \bibitem{hitchin:1987 integrable systems}
 {N. J. Hitchin}, \emph{Stable bundles and integrable systems}, Duke Math. J. \textbf{54} (1987), 91--114.

\bibitem{hitchin:1987}
{N. J. Hitchin}, \emph{The self-duality equations on a Riemann
surface}, Proc. London Math. Soc. (3) \textbf{55} (1987), 59--126.


\bibitem{lange:1983}
{H. Lange}, \textit{Universal families of extensions}, J. Algebra \textbf{83} (1983), 101--112.

\bibitem{mukai:2003}
{S. Mukai}, \emph{An Introduction to Invariants and Moduli},
Cambridge studies in advanced mathematics \textbf{81}, Cambridge University Press, 2003.

\bibitem{mundet:2000}
{I. Mundet i Riera}, \emph{A Hitchin--Kobayashi correspondence for K\"ahler
fibrations}, J. Reine Angew. Math. \textbf {528} (2000), 41--80.


\bibitem{oliveira:2008}
{A. G. Oliveira}, \textit{Higgs bundles, quadratic pairs and the topology of moduli spaces}, Ph.D. Thesis, Departamento de Matemática Pura, Faculdade de Ciências, Universidade do Porto, 2008.

\bibitem{oliveira:2010}
{A. G. Oliveira}, \textit{Representations of surface groups in the projective general linear group}, Internat. J. Math. \textbf {22} (2011), 223--279.



\bibitem{schmitt:2004}
{A. H. W. Schmitt}, \emph{A universal construction for moduli spaces of decorated vector bundles over curves}, Transform. Groups \textbf{9} (2004), 167--209.

\bibitem{schmitt:2005}
{A. H. W. Schmitt}, \emph{Moduli for decorated tuples for sheaves and representation spaces for quivers}, Proc. Indian Acad. Sci. Math. Sci. \textbf{115} (2005), 15--49.

\bibitem{schmitt:2008}
{A. H. W. Schmitt}, \emph{Geometric Invariant Theory and Decorated Principal Bundles},
Zurich Lectures in Advanced Mathematics, European Mathematical Society, 2008.

\bibitem{simpson:1992}
{C. T. Simpson}, \emph{Higgs bundles and local systems}, Inst.
Hautes \'{E}tudes Sci. Publ. Math. \textbf{75} (1992), 5--95.



\bibitem{thaddeus:1994}
{M. Thaddeus}, \emph{Stable pairs, linear systems and the Verlinde formula}, 
Invent. Math. \textbf{117} (1994), 317--353.
\end{thebibliography}

\end{document}